\documentclass[12pt]{article}

\usepackage{dsfont}
\usepackage{amsfonts,amsmath,amsthm}
\usepackage{algorithm}
\usepackage{algpseudocode}
\usepackage{color}
\usepackage{graphicx}
\usepackage{stmaryrd}        
\usepackage{cleveref} 
\usepackage{subfig}

\usepackage[paper=a4paper,dvips,top=2cm,left=2cm,right=2cm,
    foot=1cm,bottom=4cm]{geometry}

\usepackage{multicol} 
\usepackage{multirow} 

\begin{document}
\large

\title{A Power Method for Computing the Dominant  Eigenvalue of a Dual Quaternion Hermitian Matrix}
\author{ 
     Chunfeng Cui\footnote{LMIB of the Ministry of Education, School of Mathematical Sciences, Beihang University, Beijing 100191 China.
    ({\tt chungfengcui@buaa.edu.cn}).
     This author's work was supported by  Natural Science Foundation of China (Nos. 12126608, 12131004).} 
      \and \
     Liqun Qi\footnote{Department of Mathematics, School of Science, Hangzhou Dianzi University, Hangzhou 310018 China; Department of Applied Mathematics, The Hong Kong Polytechnic University, Hung Hom, Kowloon, Hong Kong
({\tt maqilq@polyu.edu.hk}).}
}
\date{\today}
\maketitle

\begin{abstract} 
In this paper, we   first study the projections onto the set of unit dual quaternions, and the set of dual quaternion vectors with unit norms. Then we propose a power method for computing the dominant eigenvalue of a dual quaternion Hermitian matrix. For {a strict dominant eigenvalue}, we show the sequence generated by the power method converges to {the} dominant eigenvalue and its corresponding eigenvector linearly. For {a general dominant eigenvalue},  we show  the standard part  of the sequence generated by the power method converges to the standard part of {the} dominant eigenvalue and its corresponding eigenvector linearly.  Based upon these, 
we reformulate the simultaneous localization and mapping (SLAM) problem as a rank-one dual quaternion completion problem. A two-block  coordinate descent method is proposed to solve this problem. 
One block     has a closed-form solution and the other block is  the best rank-one approximation problem of a dual quaternion Hermitian matrix, which can be computed by the power method.  
Numerical experiments are presented to show the efficiency of our proposed power method.

\medskip



\end{abstract}

\renewcommand{\Re}{\mathbb{R}}
\newcommand{\rank}{\mathrm{rank}}
\newcommand{\X}{\mathcal{X}}
\newcommand{\A}{\mathcal{A}}
\newcommand{\I}{\mathcal{I}}
\newcommand{\st}{\mathrm{st}}
\newcommand{\B}{\mathcal{B}}
\newcommand{\C}{\mathcal{C}}
\newcommand{\OO}{\mathcal{O}}
\newcommand{\e}{\mathbf{e}}
\newcommand{\0}{\mathbf{0}} 
\newcommand{\dd}{\mathbf{d}}
\newcommand{\ii}{\mathbf{i}}
\newcommand{\jj}{\mathbf{j}}
\newcommand{\kk}{\mathbf{k}}
\newcommand{\va}{\mathbf{a}}
\newcommand{\vb}{\mathbf{b}}
\newcommand{\vc}{\mathbf{c}}
\newcommand{\vq}{\mathbf{q}}
\newcommand{\vg}{\mathbf{g}}
\newcommand{\pr}{\vec{r}}
\newcommand{\ps}{\vec{s}}
\newcommand{\pt}{\vec{t}}
\newcommand{\pu}{\vec{u}}
\newcommand{\pv}{\vec{v}}
\newcommand{\pw}{\vec{w}}
\newcommand{\pp}{\vec{p}}
\newcommand{\pq}{\vec{q}}
\newcommand{\pl}{\vec{l}}
\newcommand{\vt}{\rm{vec}}
\newcommand{\vx}{\mathbf{x}}
\newcommand{\vy}{\mathbf{y}}
\newcommand{\vu}{\mathbf{u}}
\newcommand{\vv}{\mathbf{v}}
\newcommand{\y}{\mathbf{y}}
\newcommand{\vz}{\mathbf{z}}
\newcommand{\T}{\top}

\newcommand{\qu}[1]{\tilde{#1}}
\newcommand{\vqu}[1]{\tilde{\mathbf{#1}}}
\newcommand{\dq}[1]{\hat{#1}}
\newcommand{\vdq}[1]{\hat{\mathbf{#1}}}
\newcommand{\dqset}{\hat{\mathbb Q}}

\newtheorem{Thm}{Theorem}[section]
\newtheorem{Def}[Thm]{Definition}
\newtheorem{Ass}[Thm]{Assumption}
\newtheorem{Lem}[Thm]{Lemma}
\newtheorem{Prop}[Thm]{Proposition}
\newtheorem{Cor}[Thm]{Corollary}
\newtheorem{example}[Thm]{Example}
\newtheorem{remark}[Thm]{Remark}

\section{Introduction}\label{sec:intro}
Dual quaternion numbers and  dual quaternion matrices are  important in   robotic research, i.e., the hand-eye calibration problem \cite{Da99}, 
the simultaneous localization and mapping (SLAM) problem \cite{BS07, BLH19, CCCLSNRL16, CTDD15, WJZ13, WND00}, and the kinematic modeling and control \cite{QWL22}. 
In \cite{QL22}, Qi and Luo {studied} right and left eigenvalues {of} square dual quaternion 
matrices. If a right eigenvalue is a dual number, then it is also a left eigenvalue. In this case, this dual number is called an eigenvalue of that dual quaternion matrix. They showed that the right eigenvalues of a dual quaternion Hermitian matrix are dual numbers. Thus, they are eigenvalues. An $n$-by-$n$ dual quaternion Hermitian matrix  was shown to have exactly $n$ eigenvalues. It is positive semidefinite, or positive definite, if and only if all of its eigenvalues are nonnegative, or positive and appreciable, dual numbers, respectively. A unitary decomposition for a dual quaternion Hermitian matrix was also proposed. In \cite{QWL22}, it was shown that the eigenvalue theory of dual quaternion Hermitian matrices plays an important role in the multi-agent formation control.   However, the numerical methods for computing the eigenvalues of a dual quaternion Hermitian matrix is blank.

The power method is one of the state of the art numerical approaches for computing eigenvalues, such as  matrix eigenvalues \cite{GVL13}, matrix sparse eigenvalues \cite{YZ13}, 
tensor eigenvalues \cite{KM11}, nonnegative tensor eigenvalues \cite{NQZ10}, and quaternion matrix eigenvalues \cite{LWZZ19}.  In this paper,  we propose a power method for  computing the dominant eigenvalue of a dual quaternion Hermitian matrix. 
We first study the projections onto the set of unit dual quaternions  and the set of dual quaternion vectors with unit norms.  Then we propose a power method for computing the dual quaternion Hermitian matrix eigenvalues. 
This fills the blank of dual quaternion matrix computation.  We also define the convergence of a dual number sequence and study the convergence properties of the power method we proposed.


Dual quaternion is a powerful tool for solving the SLAM problem.  The SLAM problem aims to build a map of the environment and to simultaneously localize within  this map, which is an essential skill for mobile robots navigating in unknown environments in absence of external referencing systems such as GPS. An intuitive way to solve the SLAM problem is via graph-based formulation \cite{GKSB10}. One line of work  reformulates SLAM as a nonlinear least square and solves it by the  Gauss-Newton method \cite{KGS11}. However, as this problem is nonconvex and  Gauss-Newton may be trapped in local minimal points, a good initialization  is important \cite{CTDD15}. 
Several methods have been proposed to address the issue of global convergence, including the majorization minimization \cite{FM20} and the Lagrangian duality \cite{CRCLD15}. In this paper, we reformulate SLAM as a rank-one dual quaternion completion problem and propose a two-block  coordinate descent method to solve it. 
One block subproblem    has a closed form solution and the other block subproblem can be obtained by the dominant eigenpair. 
This connects the dual quaternion matrix theory with the SLAM problem.

The distribution of the remainder of this paper is as follows.  In the next section, we review some basic properties of  dual quaternions and dual quaternion matrices.  We also define the convergence of a dual number sequence there. 
In Section \ref{sec:proj_udq}, we show the projections onto the set of  unit dual quaternions,  and the projection onto the set of dual quaternion vectors with unit norms, can be obtained by normalization for appreciable numbers and vectors, respectively.  
In Section \ref{sec:powermethod}, we show that the sequence generated by the power method {converges} linearly to {a strict} dominant eigenvalue and eigenvector of a dual quaternion Hermitian matrix. For  {a general dominant eigenvalue}, we show   the convergence and the linear convergence rate of  the standard part of the sequence to the standard part of the dominant eigenvalue and eigenvector. 
   In Section \ref{sec:numerical}, we present numerical experiment results on computing the eigenvalues of Laplacian matrices.     
   In Section \ref{sec:SLAM}, we reformulate SLAM as a rank-one approximation model   and present a block coordinate descent method for solving it. 
   Some final remarks are made in Section \ref{sec:finial}.

\section{Dual Quaternions and Dual Quaternion Matrices}\label{sec:prelim}

\subsection{Dual quaternions}
The sets of real numbers, dual numbers, quaternion, unit quaternion, dual quaternion, and unit dual quaternion  are denoted as $\Re$, $\mathbb D$, $\mathbb Q$, $\mathbb U$, $\hat{\mathbb Q}$, and $\hat{\mathbb U}$, respectively.   Denote $\Re^3$ as $\mathbb V$.  

A {\bf quaternion} $\tilde q = [q_0, q_1, q_2, q_3]$ is a real four-dimensional vector.   We use a tilde symbol to distinguish a quaternion.   We may also write $\tilde q = [q_0, \pq]$, where $\pq = [q_1, q_2, q_3] \in \mathbb V$.  See \cite{CKJC16, Da99, WHYZ12}.  If $q_0 = 0$, then $\tilde q = [0, \pq]$ is called a {\bf vector quaternion}.
 Suppose that we have two quaternions
$\tilde p = [p_0, \pp], \tilde q = [q_0, \pq] \in \mathbb Q$, where $p = [p_1, p_2, p_3], q = [q_1, q_2, q_3] \in \mathbb V$.  The addition of $\tilde p$ and $\tilde q$ is defined as
$\tilde p + \tilde q = [p_0+q_0, \pp+\pq].$
Denote the zero element of $\mathbb Q$ as $\tilde 0 := [0, 0, 0, 0] \in \mathbb Q$.
The multiplication of $\tilde p$ and $\tilde q$ is defined by
$$\tilde p\tilde q = [p_0q_0-\pp \cdot \pq, p_0\pq+q_0\pp+\pp \times \pq],$$
where $\pp \cdot \pq$ is the dot product, i.e., the inner product of $\pp$ and $\pq$, with
$$\pp \cdot \pq \equiv \pp^\top \pq = p_1q_1+p_2q_2+p_3q_3,$$
and $\pp \times \pq$ is the cross product of $\pp$ and $\pq$, with
$$\pp \times \pq = [p_2q_3-p_3q_2,-p_1q_3+p_3q_1, p_1q_2-p_2q_1] = - \pq \times \pp.$$
Thus, in general, $\tilde p \tilde q \not = \tilde q \tilde p$, and we have $\tilde p\tilde q = \tilde q\tilde p$ if and only if $\pp \times \pq = \vec{0}$, i.e., either $\pp = \vec{0}$ or $\pq = \vec{0}$, or $\pp = \alpha \pq$ for some real number $\alpha$.

The conjugate of a quaternion $\tilde q = [q_0, q_1, q_2, q_3] \in \mathbb Q$ is defined as $\tilde q^* = [q_0, -q_1, -q_2, -q_3]$.    Let $\tilde 1 := [1, 0, 0, 0] \in \mathbb Q$.  Then for any $\tilde q \in \mathbb Q$, we have $\tilde q\tilde 1 = \tilde 1\tilde q = \tilde q$, i.e., $\tilde 1$ is the idenity element of $\mathbb Q$.     For any $\tilde p, \tilde q \in \mathbb Q$, we have
$(\tilde p\tilde q)^* = \tilde q^*\tilde p^*.$
Suppose that $\tilde p, \tilde q \in \mathbb Q$, $\tilde p \tilde q = \tilde q \tilde p = \tilde 1$.  Then we say that $\tilde p$ is invertible and its inverse is $\tilde p^{-1} = \tilde q$.

For a quaternion $\tilde q = [q_0, q_1, q_2, q_3]\in \mathbb Q$, its magnitude is defined by $$|\tilde q| = \sqrt{q_0^2 + q_1^2 + q_2^2 +q_3^2}.$$  If $|\tilde q|=1$, then it is called a {\bf unit quaternion}.   
If $\tilde p, \tilde q \in \mathbb U$, then $\tilde p\tilde q \in \mathbb U$.  For any $\tilde q \in \mathbb U$, we have
$\tilde q\tilde q^* = \tilde q^* \tilde q = \tilde 1,$
i.e., $\tilde q$ is invertible and $\tilde q^{-1} = \tilde q^*$.


A {\bf dual number} $a = a_{\st}+a_{\I}\epsilon$ consists of the standard part $a_{\st} \in \mathbb R$ and the dual part $a_{\I} \in \mathbb R$. The symbol $\epsilon$ is the infinitesimal unit.  It satisfies $\epsilon\neq 0$ and $\epsilon^2 = 0$.  
If $a_{\st}\neq0$, then $a$ is appreciable, otherwise, $a$ is infinitesimal \cite{QLY22}.

For $a = a_{\st}+a_{\I}\epsilon, b = b_{\st}+b_{\I}\epsilon \in \mathbb D$, the addition of $a$ and $b$ is defined as
$$a+b = a_{\st}+b_{\st}+(a_{\I}+b_{\I})\epsilon,$$
  the multiplication of $a$ and $b$ is defined as
$$ab = a_{\st}b_{\st}+(a_{\st}b_{\I}+a_{\I}b_{\st})\epsilon,$$ 
 and the division of $a$ and $b$, when $a_{\st}\neq 0$, or $a_{\st}= 0$ and $b_{\st}=0$ is defined as
     \begin{equation*}
         \frac{b_{\st}+b_{\I}\epsilon}{a_{\st}+a_{\I}\epsilon} =
         \left\{
         \begin{array}{cl}
             \frac{b_{\st}}{a_{\st}} + \left(\frac{b_{\I}}{a_{\st}}-\frac{b_{\st}}{a_{\st}}\frac{a_{\I}}{a_{\st}}\right)\epsilon,  & \text{ if } a_{\st}\neq 0, \\
             \frac{b_{\I}}{a_{\I}}+c \epsilon, & \ \text{if } a_{\st}= 0, b_{\st}=  0,\\
         \end{array} 
         \right.
     \end{equation*}
     where $c\in\mathbb R$ is an arbitrary  real number.   One may show that the division of dual numbers is the inverse operation of the multiplication of dual numbers.

The zero element of $\mathbb D$ is $0_D:=  0+0\epsilon$.  
The identity element of $\mathbb D$ is $1_D = 1+0\epsilon$.  
{The absolute value of $a\in\mathbb D$ is defined in \cite{QLY22} as 
\begin{equation*}
    |a|=\left\{
         \begin{array}{cl}
            |a_{\st}|+\text{sgn}(a_{\st})a_{\I} \epsilon, & \text{ if } a_{\st}\neq 0,\\
            |a_{\I}|\epsilon, & \text{otherwise}.
         \end{array} 
         \right.
\end{equation*}
   When $a\neq 0_D$, we have   $\frac{a}{|a|}=\text{sgn}(a_{\st})$ if $ a_{\st}\neq 0$, and $\frac{a}{|a|}=\text{sgn}(a_{\I})+c\epsilon$ for any $c\in\mathbb R$ if $ a_{\st}=0$.}

The total order of   dual numbers was defined in \cite{QLY22}: we say that $a>b$ if 
$a_{\st}>b_{\st}$, or $a_{\st}=b_{\st}$ and $a_{\I}>b_{\I}$. 
Let $a>0_D$ be a {nonnegative} dual number. Then the square root of $a$ is defined by \cite{QLY22}
\begin{equation*}
    \sqrt{a}=\sqrt{a_{\st}}+\frac{a_{\I}}{2\sqrt{a_{\st}}}\epsilon.
\end{equation*}


 Let $\{a_k = a_{k, \st}+ a_{k, \I}\epsilon : k = 1, 2, \cdots \}$ be a dual number sequence.  We say that this sequence is convergent and has a limit $a = a_{\st} + a_{\I}\epsilon$ if both its standard part sequence $\{ a_{k, st} : k = 1, 2, \cdots \}$ and its dual part sequence $\{ a_{k, \I} : k = 1, 2, \cdots \}$ are convergent and have limits $a_{\st}$ and $a_{\I}$ respectively.

Given a dual number sequence $\{ a_k : k = 1, 2, \cdots \}  \subset \mathbb D$ and a real number sequence $\{ c_k : k = 1, 2, \cdots \}  \subset \mathbb R$, we introduce a big $O_D$ notation as follows,  
\begin{equation}\label{equ:O_D}
    a_k = O_D(c_k) \text{ if } a_{k,\st} = O(c_k) \text{ and }a_{k,\I} = O(c_k). 
\end{equation}
Clearly, in this case, if the real number sequence $\{c_k : k = 1, 2, \cdots \}$ is convergent, then the dual number sequence $\{ a_k : k = 1, 2, \cdots \}$ is also convergent. Furthermore, if $c_k=O(c^k h(k))$   for a real positive number $c<1$ and a polynomial $h(k)$, then we denote $c_k=\tilde{O}(c^k)$. 
Similarly,  we denote 
\begin{equation}\label{equ:titleO_D}
    a_k = \tilde O_D(c_k) \text{ if } a_{k,\st} = \tilde O(c_k) \text{ and }a_{k,\I} = \tilde O(c_k). 
\end{equation}

A {\bf dual quaternion} $\hat q =  \tilde q_{\st}+ \tilde q_{\I}\epsilon$ consists of two quaternions $\tilde q_{\st}$, the standard part of $\hat q$, and $\tilde q_{\I}
$, the dual part of $\hat q$.      We use a hat symbol to distinguish a dual quaternion.  
If $\qu q_{\st}\neq \qu0$, then $\dq q$ is appreciable, otherwise, $\dq q$ is infinitesimal \cite{QLY22}.

Let $\hat p =  \tilde p_{\st}+\tilde p_{\I} \epsilon, \hat q =  \tilde q_{\st}+\tilde q_{\I}\epsilon \in \hat {\mathbb Q}$.   Then the sum of $\hat p$ and $\hat q$ is
$$\hat p + \hat q =  (\tilde p_{\st} + \tilde q_{\st})+(\tilde p_{\I}+\tilde q_{\I})\epsilon$$
and the product of $\hat p$ and $\hat q$ is
$$\hat p\hat q =  \tilde p_{\st}\tilde q_{\st}+(\tilde p_{\st}\tilde q_{\I} + \tilde p_{\I}\tilde q_{\st})\epsilon.$$
Again, in general, $\hat p \hat q \not = \hat q\hat p$.   The conjugate of $\hat p =  \tilde p_{\st}+ \tilde p_{\I}\epsilon$ is $\hat p^* =  \tilde p_{\st}^*+ \tilde p_{\I}^*\epsilon$. {Denote $\hat 0 = \tilde 0 + \tilde 0\epsilon$ and $\hat 1 = \tilde 1 + \tilde 0\epsilon$.}

A dual quaternion $\hat q =  \tilde q_{\st}+\tilde q_{\I}\epsilon$ is called a {\bf unit dual quaternion} if 
\begin{equation*} 
|\tilde q_{\st}|=1 \text{ and } \tilde q_{\st}\tilde q_{\I}^* + \tilde q_{\I}\tilde q_{\st}^* = \tilde 0.
\end{equation*}
The magnitude of a dual quaternion $\hat q\in \hat{\mathbb Q}$ is \cite{QLY22}
\begin{equation*}
    |\hat q|=\left\{
    \begin{array}{cc}
       |\tilde q_{\st}|+\frac{sc(\tilde q_{\st}^*\tilde q_{\I})}{|\tilde q_{\st}|} \epsilon,  &  \text{ if } \tilde q_{\st}\neq \qu 0,\\
       |\tilde q_{\I}|\epsilon,  & \text{otherwise}.
    \end{array}
    \right.
\end{equation*}
Here,  $sc(\tilde q)=\frac12 (\tilde q+\tilde q^*)$ is the scalar part of $\tilde q$. 

Further, the $\ell_{2^*}$-norm similar to \cite{WDW23} can be defined by 
\begin{equation*}
    |\hat q|_*=\left\{
    \begin{array}{cc}
       |\tilde q_{\st}|+\frac{|\tilde q_{\I}|^2}{2|\tilde q_{\st}|} \epsilon,  &  \text{ if } \tilde q_{\st}\neq 0,\\
       |\tilde q_{\I}|\epsilon,  & \text{otherwise}.
    \end{array}
    \right.
\end{equation*}

For any  dual quaternion number $\dq q=\qu q_{\st}+\qu q_{\I}\epsilon$ and dual number  $a=a_{\st}+a_{\I}\epsilon$ with $a_{\st}\neq 0$, or $a_{\st}= 0$ and $\qu q_{\st}=\qu 0$, there is 
     \begin{equation*}
         \frac{\qu q_{\st}+\qu q_{\I}\epsilon}{a_{\st}+a_{\I}\epsilon} =
         \left\{
         \begin{array}{cl}
             \frac{\qu q_{\st}}{a_{\st}} + \left(\frac{\qu q_{\I}}{a_{\st}}-\frac{\qu q_{\st}}{a_{\st}}\frac{a_{\I}}{a_{\st}}\right)\epsilon,  & \text{ if } a_{\st}\neq 0, \\
             \frac{\qu q_{\I}}{a_{\I}}+\tilde c \epsilon, & \ \text{if } a_{\st}= 0, \qu q_{\st}=\qu 0,\\
         \end{array} 
         \right.
     \end{equation*}
     where $\tilde c\in\mathbb Q$ is an arbitrary quaternion number. 
     
\subsection{Dual quaternion matrices}

We denote the set of $m$-by-$n$ real,  quaternion, unit quaternion,  dual quaternion, and unit dual quaternion matrices as $\mathbb R^{m\times n}$, $\mathbb Q^{m\times n}$, $\mathbb U^{m\times n}$, $\hat{\mathbb Q}^{m\times n}$, and $\hat{\mathbb U}^{m\times n}$, respectively. 
We use $\vqu{O}_{m\times n}$ and $\vdq O_{m\times n}$ to denote the $m$-by-$n$ zero     quaternion  and zero dual quaternion matrices, and $\vqu{I}_{m\times m}$ and $\vdq I_{m\times m}$ to denote the $m$-by-$m$ identity      quaternion  and identity dual quaternion matrices, respectively.
{Given $\vdq Q\in\dqset^{m\times n}$. If $\vqu Q_{\st}\neq \vqu O_{m\times n}$, then $\vdq Q$ is appreciable, otherwise, $\vdq Q$ is infinitesimal.}
 
The $2$-norm of a dual quaternion vector $\hat {\mathbf x}=(\hat x_{i})\in \hat{\mathbb Q}^{n\times 1}$ is {
\begin{equation*}\label{equ:dq2norm}
    \|\hat {\mathbf x}\|_2=\left\{\begin{array}{cl}
       \sqrt{{\sum_{i=1}^n} |\hat x_{i}|^2},   & \text{ if }\vqu x_{\st}\neq \vqu O_{n\times 1}, \\
       \sqrt{{\sum_{i=1}^n} |\qu x_{i,\I}|^2} \epsilon,  &  \text{ if } \vqu x_{\st}=\vqu O_{n\times 1} \text{ and }\vdq x  = \vqu x_{\I}\epsilon.
    \end{array} 
    \right.
\end{equation*}
}
A vector $\hat {\mathbf x}=(\hat x_{i}) \in\hat{\mathbb U}^{n\times 1}$ is called a unit dual quaternion vector if  each element in  $\hat {\mathbf x}$ is a unit dual quaternion number. The set of $n \times 1$ quaternion vectors with unit $2$-norms is denoted by $\mathbb Q^{n \times 1}_2$. 
A matrix $\hat {\mathbf Q}=(\hat q_{i,j}) \in\hat{\mathbb U}^{m\times n}$ is called a unit dual quaternion matrix if  each element in $\hat {\mathbf Q}$ is a unit dual quaternion number.

Given a matrix $\hat {\mathbf Q}=(\hat{q}_{i,j}) \in\dqset^{m\times n}$, the transpose of $\hat {\mathbf Q}$ is $\hat {\mathbf Q}^\top=(\hat{q}_{j,i}) \in{\hat{\mathbb Q}}^{n\times m}$, and the conjugate transpose of $\hat {\mathbf Q}$ is $\hat {\mathbf Q}^*=(\hat{q}^*_{j,i}) \in\hat{\mathbb Q}^{n\times m}$. 
If $\hat {\mathbf Q}=\hat {\mathbf Q}^*$, then it is a dual quaternion Hermitian matrix, and both its standard part and infinitesimal part are quaternion Hermitian matrices.

The $F$-norm of a dual quaternion matrix $\hat {\mathbf Q}=(\hat Q_{ij})\in \hat{\mathbb Q}^{m\times n}$ is 
{
\begin{equation*}\label{equ:dqFnorm}
    \|\hat {\mathbf Q}\|_F=\left\{\begin{array}{cl}
       \|\vqu Q_{\st}\|_F+\frac{sc(tr(\vqu Q_{\st}^*\vqu Q_{\I}))}{\|\vqu Q_{\st}\|_F}\epsilon,   &  \text{ if }\vqu Q_{\st}\neq \vqu O_{m\times n}, \\
       \|\vqu Q_{\I}\|_F\epsilon,   &  \text{ if }\vqu Q_{\st}=\vqu O_{m\times n} \text{ and } \vdq Q=\vqu Q_{\I}\epsilon,
    \end{array}\right.  
\end{equation*} 
}
and the $F^*$-norm  is 
\begin{equation*}\label{equ:dqF*norm}
    \|\hat {\mathbf Q}\|_{F^*}=\left\{\begin{array}{cl}
     \|\vqu Q_{\st}\|_F+\frac{\|\vqu Q_{\I}\|^2}{2\|\vqu Q_{\st}\|_F}\epsilon,   &  \text{ if }\vqu Q_{\st}\neq \vqu O_{m\times n}, \\
       \|\vqu Q_{\I}\|_F\epsilon,   &  \text{ if }\vqu Q_{\st}=\vqu O_{m\times n} \text{ and } \vdq Q=\vqu Q_{\I}\epsilon,
    \end{array}\right.
\end{equation*}

For convenience in the numerical experiments, we define the $2^R$-norm of a dual quaternion vector $\hat {\mathbf x}=(\hat x_{i}) \in\hat{\mathbb Q}^{n\times 1}$  as 
\begin{equation*}
    \|\hat {\mathbf x}\|_{2^R} = \sqrt{\|\vqu x_{\st}\|_2^2 + \|\vqu x_{\I}\|_2^2}, 
\end{equation*}
and the $F^R$-norm of a dual quaternion matrix $\vdq Q=(\dq q_{ij}) \in\hat{\mathbb Q}^{m\times n}$  as 
\begin{equation*}
    \|\hat {\mathbf Q}\|_{F^R} = \sqrt{\|\vqu Q_{\st}\|_F^2 + \|\vqu Q_{\I}\|_F^2}, 
\end{equation*}
respectively.

\bigskip

\section{Projection onto the Set of Dual Quaternions with Unit Norms}\label{sec:proj_udq}


Given a quaternion number $\qu q\neq \qu 0$ and a quaternion vector $\vqu q\neq \vqu O_{n\times 1}$, we have 
\begin{equation*}
    P_{\mathbb U}(\qu q) = \frac{\qu q}{|\qu q|} \text{ and } P_{\mathbb Q^{n \times 1}_2}(\vqu q) = \frac{\vqu q}{\|\vqu q\|_2}. 
\end{equation*}
Here, $\mathbb Q^{n \times 1}_2$ is the set of $n$ dimensional quaternion vectors with unit $2$-norms. 
In the following, we will show that for the dual quaternion number and the dual quaternion vector, the projection has a similar formulation. 


\begin{Thm} \label{prop:proj_qu} 
Given a dual quaternion number $\dq q=\qu q_{\st}+\qu q_{\I}\epsilon\in\dqset$, we have the following conclusions.   
    \begin{itemize}
        \item[(i)]  If $\qu q_{\st}\neq \qu 0$, the normalization of $\dq q$ is the projection of $\dq q$ onto the unit dual quaternion set. Namely,   
        \begin{equation}\label{equ:dq_normalize}
         \dq u = \frac{\dq q}{|\dq q|} \text{ with } \qu u_{\st}=\frac{\qu q_{\st}}{|\qu q_{\st}|},\ \qu u_{\I}= \frac{\qu q_{\I}}{|\qu q_{\st}|}-\frac{\qu q_{\st}}{|\qu q_{\st}|}sc\left(\frac{\qu q_{\st}^*}{|\qu q_{\st}|}\frac{\qu q_{\I}}{|\qu q_{\st}|}\right)  
    \end{equation}
    is a unit dual quaternion number and     \begin{equation}\label{equ:prob_dq_proj}
       \dq u\in\mathop{\arg\min}_{\dq v\in\hat{\mathbb U}}\quad |\dq v-\dq q|^2.
    \end{equation}
    \item[(ii)]  If $\qu q_{\st}= \qu 0$ and $\qu q_{\I}\neq \qu 0$,    then the projection of $\dq q$ onto the unit dual quaternion set is  $\dq u = \qu u_{\st} + \qu u_{\I} \epsilon$, where
        \begin{equation}\label{equ:dq_st0_normalize}
           \qu u_{\st}=\frac{\qu q_{\I}}{|\qu q_{\I}|}, \ \qu u_{\I} \text{ is any quaternion number  satisfying } {sc(\qu q_{\I}^*\qu u_{\I})=0}. 
    \end{equation}
    \end{itemize}
\end{Thm}
\begin{proof}
    $(i)$ The normalization formula in equation \eqref{equ:dq_normalize} follows from 
      \begin{equation*}
        \frac{\dq q}{|\dq q|} = \frac{\qu q_{\st}+\qu q_{\I}\epsilon}{|\qu q_{\st}| + \frac{sc(\qu q^*_{\st}\qu q_{\I})}{|\qu q_{\st}|}\epsilon}=\frac{\qu q_{\st}}{|\qu q_{\st}|}+\left(\frac{\qu q_{\I}}{|\qu q_{\st}|}-\frac{\qu q_{\st}}{|\qu q_{\st}|}sc\left(\frac{\qu q_{\st}^*}{|\qu q_{\st}|}\frac{\qu q_{\I}}{|\qu q_{\st}|}\right)\right)\epsilon.  
    \end{equation*}
    By direct computations, we have $\qu u_{\st}^*\qu u_{\st}=\qu 1$ and $\qu u_{\st}^*\qu u_{\I}+\qu u_{\I}^*\qu u_{\st}=0$. Furthermore, by
     $$|\dq v-\dq q|^2 = |\qu v_{\st}-\qu q_{\st}|^2 + 2sc((\qu v_{\st}-\qu q_{\st})^*(\qu v_{\I}-\qu q_{\I}))\epsilon$$ and the   definition of the  total order of dual numbers, we have 
     \begin{equation}\label{equ:proj_stand}
         \tilde u_{\st} = \mathop{\arg\min}\limits_{\qu v_{\st}\in\mathbb U} \ |\qu v_{\st}-\qu q_{\st}|^2= \frac{\qu q_{\st}}{|\qu q_{\st}|}.
     \end{equation}
    Moreover, we have $sc((\qu u_{\st}-\qu q_{\st})^*(\qu v_{\I}-\qu q_{\I}))=-sc((\qu u_{\st}-\qu q_{\st})^*\qu q_{\I})$
     for all $\qu v_{\I}$ satisfying $sc(\qu q_{\st}^*\qu v_{\I})=0$. In other words, {any $\dq u$  with the standard part  $\qu u_{\st}=\frac{\qu q_{\st}}{|\qu q_{\st}|}$  and the infinitesimal part   satisfying $sc(\qu q_{\st}^*\qu u_{\I})=0$ is an optimal solution.}  Hence, we conclude that $\dq u$ in \eqref{equ:dq_normalize} is an optimal solution of \eqref{equ:prob_dq_proj}. 

    $(ii)$  When $\qu q_{\st}=0$, we have 
    $$|\dq v-\dq q|^2=|\qu v_{\st}|^2+2sc(\qu v_{\st}^*(\qu v_{\I}-\qu q_{\I}))\epsilon = 1-2sc(\qu v_{\st}^*\qu q_{\I})$$
     for any $\dq v\in\hat{\mathbb U}$. 
    By the definition of the total order of dual numbers, there is  
    \begin{equation*}
         \tilde u_{\st} = \mathop{\arg\min}\limits_{\qu v_{\st}\in\mathbb U}-sc(\qu v_{\st}^*\qu q_{\I}). 
     \end{equation*}
    Hence, $\tilde u_{\st}=\frac{\qu q_{\I}}{|\qu q_{\I}|}$. 
    As the objective function is independent of $\qu v_{\I}$, the infinitesimal part of $\dq u$ can be any quaternion number satisfying $sc(\qu q_{\I}^*\qu u_{\I})=0$. 
    This completes the proof. 
\end{proof}

   For problem \eqref{equ:prob_dq_proj},  any $\qu u_{\I}$ satisfying $sc(\qu q_{\st}^*\qu u_{\I})=0$ is  an optimal solution. However, the choice of $\qu u_{\I}$ in \eqref{equ:dq_normalize}   is  geometric meaningful, as shown in the following proposition. 
    
    \begin{Prop}\label{prop:proj_dual}
        Suppose $\qu q_{\st}\neq\qu 0$. Then $\qu u_{\I}$ defined by  \eqref{equ:dq_normalize} is a solution to the following optimization problem, 
    \begin{equation}\label{equ:uIprob}
        \min_{\qu v\in \mathbb Q}\quad \left|\qu v -\frac{\qu q_{\I}}{|\qu q_{\st}|}\right|^2\ \mathrm{ s.t. }\ sc(\qu v^*\qu q_{\st})=0. 
    \end{equation}
    Furthermore, $\dq u = \frac{\dq q}{|\dq q|}$ is the   optimal solution to the following problem 
     \begin{equation*} 
       \dq u=\mathop{\arg\min}_{\dq v\in\hat{\mathbb U}}\quad |\qu v_{\st}-\qu q_{\st}|^2+\left|\qu v_{\I}-\frac{\qu q_{\I}}{|\qu q_{\st}|}\right|^2\epsilon.
    \end{equation*}
    \end{Prop}
    \begin{proof}
    {Problem \eqref{equ:uIprob} is a convex quaternion optimization problem  \cite{QLWZ22} since  the objective function    is convex and the constraint is linear.   Then it follows from  the first-order optimality conditions given in Theorem 4.3 of \cite{QLWZ22}, there is a Lagrange multiplier $\lambda\in\mathbb R$ such that 
    \begin{equation*}
        \qu v -\frac{\qu q_{\I}}{|\qu q_{\st}|}+\lambda \qu q_{\st}=0\  \text{ and } \ sc(\qu v^*\qu q_{\st})=0.
    \end{equation*}
    By multiplying $\qu q_{\st}^*$ at the both sides of the first equation, there is $\lambda = \frac{sc(\qu q_{\st}^*\qu q_{\I})}{|\qu q_{\st}|^3}$. Consequently, we have 
    \begin{equation*}
        \qu v =\frac{\qu q_{\I}}{|\qu q_{\st}|}-\lambda \qu q_{\st}=\frac{\qu q_{\I}}{|\qu q_{\st}|}-   \frac{sc(\qu q_{\st}^*\qu q_{\I})}{|\qu q_{\st}|^3} \qu q_{\st}.
    \end{equation*}
    This completes the proof. 
    }
    \end{proof}

 In the case that $\qu q_{\st}=\qu 0$ and $\qu q_{\I}\neq\qu 0$, the normalization of $\dq q$ is 
 \begin{equation*}
     \frac{\dq q}{|\dq q|} = \frac{\qu q_{\I}}{|\qu q_{\I}|}+\qu v_{\I}\epsilon,
 \end{equation*}
 where $\qu v_{\I}$ can be any quaternion number. 
 The standard part is the same as that of \eqref{equ:dq_st0_normalize} and the dual part is different because  \eqref{equ:dq_st0_normalize} needs an additional condition $sc(\qu q_{\I}^*\qu u_{\I})=0$.

Similarly, we can show the projection onto the set of dual quaternion vectors with unit norms also have  closed-form solutions. 
\begin{Thm}\label{prop:proj_quvector}
Given a dual quaternion vector $\vdq q\in\dqset^{n\times 1}$, we have 
 \begin{itemize}
        \item[(i)]  If $\vqu q_{\st}\neq \vqu O_{n\times 1}$, the normalization of $\vdq q$ is the projection of $\vdq q$ onto the set of dual quaternion vectors with unit norms. Namely,   
        \begin{equation}\label{equ:vdq_normalize}
         \vdq u=\frac{\vdq q}{\|\vdq q\|_2} \text{ with } \vqu u_{\st}=\frac{\vqu q_{\st}}{\|\vqu q_{\st}\|_2},\ \vqu u_{\I}=\frac{\vqu q_{\I}}{\|\vqu q_{\st}\|_2}-\frac{\vqu q_{\st}}{\|\vqu q_{\st}\|_2}sc\left(\frac{\vqu q_{\st}^*}{\|\vqu q_{\st}\|_2}\frac{\vqu q_{\I}}{\|\vqu q_{\st}\|_2}\right) 
    \end{equation}
    is a   dual quaternion vector with  unit norm and 
    \begin{equation*}\label{equ:prob_vdq_proj}
       \vdq u \in \mathop{\arg\min}_{\vdq v\in{\hat{\mathbb Q}}_2^{n\times 1}}\quad \|\vdq v-\vdq q\|_2^2.
    \end{equation*}
    \item[(ii)]  If $\vqu q_{\st}= \vqu O_{n\times 1}$,  the   projection of $\vdq q$ onto the set of dual quaternion vectors with unit norms is   $\vdq u=\vqu u_{\st}+\vqu u_{\I}\epsilon$ satisfying 
        \begin{equation*}\label{equ:vdq_st0_normalize}
          \vqu u_{\st}=\frac{\vqu q_{\I}}{\|\vqu q_{\I}\|_2}, \ \vqu u_{\I} \text{ is any quaternion satisfying } sc(\vqu q_{\I}^*\vqu u_{\I})=0. 
    \end{equation*}
    \end{itemize}
\end{Thm}

\begin{proof}
    The proof of this theorem is similar to the proof of Theorem \ref{prop:proj_qu}.   We do not repeat it here. 
\end{proof}

{
Similar to Proposition \ref{prop:proj_dual}, we have the following result. 

 \begin{Prop}\label{prop:proj_vdual}
        If $\vqu q_{\st}\neq \vqu O_{n\times 1}$, $\vqu u_{\I}$ in \eqref{equ:vdq_normalize} is a solution to the following optimization problem, 
    \begin{equation*}\label{equ:vuIprob}
        \min_{\vqu v\in \mathbb Q_2^{n\times 1}}\quad \left\|\vqu v -\frac{\vqu q_{\I}}{\|\vqu q_{\st}\|_2}\right\|_2^2\ \mathrm{ s.t. }\ sc(\vqu v^*\vqu q_{\st})=0. 
    \end{equation*}
    Furthermore, $\vdq u = \frac{\vdq q}{\|\vdq q\|_2}$ is the   optimal solution to the following problem 
     \begin{equation*} 
       \vdq u=\mathop{\arg\min}_{\vdq v\in\hat{\mathbb Q}_2^{n\times 1}}\quad \|\vqu v_{\st}-\vqu q_{\st}\|_2^2+\left\|\vqu v_{\I}-\frac{\vqu q_{\I}}{\|\vqu q_{\st}\|_2}\right\|_2^2\epsilon.
    \end{equation*}
    \end{Prop}
    \begin{proof}
    The proof is similar to that of Proposition \ref{prop:proj_dual} and we omit it here. 
    \end{proof}
}
\bigskip
\section{The Power Method for Computing the Dominant Eigenvalue of a Dual Quaternion Hermitian Matrix}\label{sec:powermethod}

For a quaternion matrix $\vqu Q$,   the power method  can return the 
  eigenvalue with the maximum absolute value and its associated eigenvector \cite{LWZZ19}. We now study the power method for  dual quaternion  Hermitian matrices. 

\subsection{The dominant eigenvalues of dual quaternion Hermitian matrices}
Based on Theorem 4.2 in \cite{QL22}, a dual quaternion Hermitian matrix $\vdq Q\in\dqset^{n\times n}$ has exactly $n$  eigenvalues, which are  dual numbers. It also has  $n$ eigenvectors $\vdq u_1,\dots,\vdq u_n$,  which  are orthonormal vectors, i.e.,
\begin{equation*}
    \vdq Q\vdq u_i = \vdq u_i\lambda_i \text{ and }\vdq u_i^*\vdq u_j=\left\{
    \begin{array}{cc}
       \dq 1,  & \text{ if }i=j; \\
       \dq 0  &  \text{ otherwise},
    \end{array}
    \right.
    \forall\, i,j,=1,\dots,n.
\end{equation*}
Furthermore, $\{ \vdq u_i : i = 1, \dots, n \}$ forms an orthonormal basis of $\dqset^{n\times 1}$.
For any  unit dual quaternion numbers $\dq q_i\in\hat{\mathbb{U}}$, $i=1,\dots,n$, if we multiply $\dq q_i$ on the right of $\vdq u_i$, there is
\begin{equation*}
   \vdq Q\vdq u_i\dq q_i=\vdq u_i\lambda_i\dq q_i= (\vdq u_i\dq q_i)\lambda_i \text{ and }\dq q_i^*\vdq u_i^*\vdq u_j\dq q_j=\left\{
    \begin{array}{cc}
       \dq 1,  & \text{ if }i=j; \\
       \dq 0  &  \text{ otherwise}
    \end{array}
    \right.
    \forall\, i,j,=1,\dots,n, 
\end{equation*}
where the second  equality follows from that a dual number is commutative with a dual quaternion matrix. 
Hence, for a dual quaternion Hermitian matrix, 
the unit norm eigenvectors $\{\vdq u_i\}_{i=1}^n$,  which form an orthonormal basis of $\dqset^{n\times 1}$,   are not unique. 

 Suppose that we have a dual quaternion Hermitian matrix 
$\vdq Q\in\dqset^{n\times n}$.  We say that an eigenvalue $\lambda_1$ of $\vdq Q$ is a {\bf dominant eigenvalue} of $\vdq Q$, and an eigenvector corresponding to $\lambda_1$ a {\bf dominant eigenvector}, if for any eigenvalue $\lambda_j $ of $\vdq Q$, we have 
$$|\lambda_{1, st}| \ge |\lambda_{j, st}|, \ \forall\, j=2,\dots,n.$$
We say that an eigenvalue $\lambda_1$ of $\vdq Q$ is a {\bf strict dominant eigenvalue} of $\vdq Q$, with multiplicity $l$, if $\vdq Q$ has eigenvalues 
$\lambda_j$ for $j = 1, \cdots, n$, and they satisfy
\begin{equation*} 
   \lambda_{1}=\dots= \lambda_{l} \text{ and }|\lambda_{1,\st}|>|\lambda_{l+1,\st}|\ge\dots\ge|\lambda_{n,\st}|\ge0.
\end{equation*} 

{\bf Assumption A}   The standard parts of the dominant eigenvalues have the same sign.

If $\vdq Q$ has a strict dominant eigenvalue, then $\vdq Q$ satisfies Assumption A.   On the other hand, if $\vdq Q$ does not satisfy Assumption A.  Let $\vdq P = \vdq Q + \alpha \vdq I_{n\times n}$, where $\alpha$ is a nonzero real number.  Then $\dq P$ must satisfy Assumption A, and $\lambda$ is an eigenvalue of $\vdq Q$ if and only if $\lambda + \alpha$ is an eigenvalue of $\vdq P$.  Thus, it is adequate that we may assume that $\vdq Q$ satisfies Assumption A. 

\subsection{Computing the strict dominant eigenvalues by the power method}

If $\vqu Q_{\st}\neq \vqu O$, then at least one eigenvalue has 
 nonzero standard part. 
Given any initial    dual quaternion vector $\vdq v^{(0)}$ with  unit norm, the power method   computes 
\begin{equation}\label{equ:DQM_powermethod}
     \vdq{y}^{(k)} = \vdq{Q}\vdq{v}^{(k-1)}, \quad \lambda^{(k-1)}=(\vdq v^{(k-1)})^*\vdq y^{(k)},\quad \vdq{v}^{(k)} =\frac{\vdq{y}^{(k)}}{\|\vdq{y}^{(k)}\|_2}  
 \end{equation}
iteratively. 
This process is repeated until convergent {or the maximal iteration number is reached.}

\begin{algorithm}
\caption{Power method for computing the dominant eigenvalues of a dual quaternion Hermitian matrix}\label{alg:power}
\begin{algorithmic}
\Require The Hermitian matrix $\vdq{Q}$,   the initial point $\vdq v^{(0)}$, the maximal iteration number  $k_{\max}$, and the tolerance $\delta$. 
\For{$k=1,\dots, k_{\max}$}
 \State Update $\vdq y^{(k)}$   by $\vdq{y}^{(k)} = \vdq{Q}\vdq{v}^{(k-1)}$. 
 \State Update $\lambda^{(k-1)}=(\vdq v^{(k-1)})^*\vdq y^{(k)}$. 
 \If{$\|\vdq y^{(k)}-\vdq v^{(k-1)}\lambda^{(k-1)}\|_{2^R}\le \delta\times \|\vdq Q\|_{F^R}$}
 \State Stop.
 \EndIf 
 \State Update  $\vdq v^{(k)}$ by $\vdq{v}^{(k)} =\frac{\vdq{y}^{(k)}}{\|\vdq{y}^{(k)}\|_2}$. 
 \EndFor
\State \textbf{Output:} $\vdq v^{(k-1)}$   and $\lambda^{(k)}$. 
\end{algorithmic}
\end{algorithm}


Similar to the power method for the real matrices and quaternion matrices, we show $\vdq v^{(k)}$ converges to the eigenvector corresponding to  a strict dominant eigenvalue linearly. 

\begin{Thm}\label{thm:power_multiplicity_l}
 Given a dual quaternion Hermitian matrix $\vdq Q\in\dqset^{n\times n}$ and a dual quaternion vector $\vdq v^{(0)}=\sum_{j=1}^n\vdq u_j\dq \alpha_j$ with unit norm,  where $\vdq u_1,\dots,\vdq u_n$ are orthonormal eigenvectors of $\vdq Q$. Suppose $\sum_{j=1}^l|\qu \alpha_{j,\st}| \neq  0$, $\vqu Q_{\st}\neq \vqu O_{n\times n}$, and $\vdq Q$   has a strict dominant eigenvalue with multiplicity  $l$, i.e., 
\begin{equation*} 
    \lambda_{1}=\dots= \lambda_{l} \text{ and }|\lambda_{1,\st}|>|\lambda_{l+1,\st}|\ge\dots\ge|\lambda_{n,\st}|\ge0.  
\end{equation*}  
 Then  the sequence generated by the power method \eqref{equ:DQM_powermethod} satisfies 
\begin{equation}\label{equ:powermethod_vk}
    \vdq{v}^{(k)} = s^k\sum_{j=1}^l\vdq u_j \dq \gamma_j \left(1_D + \tilde O_D\left(\left|\frac{\lambda_{l+1,\st}}{\lambda_{1,\st}}\right|^k\right)\right), 
\end{equation}
where $s=sgn(\lambda_{1,\st})$, $\dq \gamma_j=\frac{\dq \alpha_j}{\sqrt{\sum_{i=1}^l |\dq\alpha_i|^2}}$, and 
\begin{equation*} 
    \lambda^{(k)}=
    \left(\vdq{v}^{(k)}\right)^*\vdq y^{(k+1)} = \lambda_1 \left(1_D + \tilde O_D\left(\left|\frac{\lambda_{l+1,\st}}{\lambda_{1,\st}}\right|^{2k}\right)\right).  
\end{equation*}
 In other words, 
the sequence  $\vdq{v}^{(k)}s^k$  converges to a strict dominant eigenvector $\sum_{j=1}^l\vdq u_j \dq \gamma_j$ at a rate of $\tilde O_D\left(\left|\frac{\lambda_{l+1,\st}}{\lambda_{1,\st}}\right|^{k}\right)$   and the sequence $\lambda^{(k)}$ converges to the strict dominant eigenvalue $\lambda_1$ at a rate of $\tilde O_D\left(\left|\frac{\lambda_{l+1,\st}}{\lambda_{1,\st}}\right|^{2k}\right)$.
 
\end{Thm} 
\begin{proof}
Let $\vdq v^{(0)}$ be a dual quaternion vector with unit norm. Since $\vdq u_1,\dots,\vdq u_n$ form a basis in $\dqset^n$, we have 
    \begin{equation*}
        \vdq v^{(0)}=\sum_{j=1}^n \vdq u_j\dq \alpha_j,
    \end{equation*}
    where $\dq \alpha_j,j=1,\dots,n$ are dual quaternion numbers, $\sum_{j=1}^n  \dq \alpha_j^*\dq \alpha_j=\dq 1$, and $\sum_{j=1}^l|\qu \alpha_{j,\st}|\neq \qu 0$.  
    Combing this with \eqref{equ:DQM_powermethod}, we have 
    \begin{equation}\label{equ:vkyk}
        \vdq v^{(k)} = \frac{\sum_{j=1}^n \vdq u_j \lambda_j^k\dq\alpha_j}{\sqrt{\sum_{j=1}^n |\lambda_j^k\dq\alpha_j|^2}}  \text{ and }\vdq y^{(k+1)} = \frac{\sum_{j=1}^n \vdq u_j \lambda_j^{k+1}\dq\alpha_j}{\sqrt{\sum_{j=1}^n |\lambda_j^k\dq\alpha_j|^2}}. 
    \end{equation}
     
     By direct derivations, we have $|\lambda_j^k\dq\alpha_j| = |\lambda_j^k||\dq\alpha_j|$. If  $\lambda_{j,\st}\neq 0$,   we have  \begin{equation*}\lambda_j^k=\lambda_{j,\st}^k+k\lambda_{j,\st}^{k-1}\lambda_{j,\I}\epsilon = \lambda_{j,\st}^k(1+k\beta_j\epsilon)
     \text{ and } |\lambda_j^k|=|\lambda_{j,\st}|^{k}(1+k\beta_j\epsilon),
     \end{equation*} 
    where      $\beta_j=\lambda_{j,\st}^{-1}\lambda_{j,\I}$. 
   If $\lambda_{j,\st}= 0$,  we have $\lambda_j^k=0_D$ for all $k\ge2$.  
When $k$ goes to infinity, we have  
    \begin{eqnarray*}
       &&  \sqrt{\sum_{j=1}^n |\lambda_j^k\dq\alpha_j|^2}-|\lambda_1|^{k}\sqrt{\sum_{j=1}^l|\dq\alpha_j|^2} \\
    &=& \frac{\sum_{j=l+1}^n |\lambda_j^k\dq\alpha_j|^2}{\sqrt{\sum_{j=1}^n |\lambda_j^k\dq\alpha_j|^2}+|\lambda_1|^{k}\sqrt{\sum_{j=1}^l|\dq\alpha_j|^2}}\\
    &=& \frac{\sum_{j=l+1}^n |\lambda_{j,\st}|^{2k}(1+2k\beta_j\epsilon)|\dq\alpha_j|^2}{ \left(\sum_{j=1}^n  |\lambda_{j,\st}|^{2k}(1+2k\beta_j\epsilon)|\dq\alpha_j|^2\right)^{1/2}+|\lambda_{1,\st}|^{k}(1+k\beta_1\epsilon)\sqrt{\sum_{j=1}^l|\dq\alpha_j|^2}}\\
    &=& O_D(|\lambda_{l+1,\st}|^k), 
    \end{eqnarray*}
     where $O_D(\cdot)$ is defined by \eqref{equ:O_D}. 
     Hence, 
 \begin{equation*}
        \frac{\lambda_j^k\dq\alpha_j}{\sqrt{\sum_{j=1}^n |\lambda_j^k\dq\alpha_j|^2}}=\frac{\lambda_j^k\dq\alpha_j}{|\lambda_1|^k\sqrt{\sum_{j=1}^l|\dq\alpha_j|^2} + O_D(|\lambda_{l+1,\st}|^k)}.   
    \end{equation*}

    If $j\in\{1,\dots,l\}$, it holds that 
\begin{eqnarray*} 
     \frac{\lambda_1^k\dq\alpha_j}{\sqrt{\sum_{j=1}^n |\lambda_j^k\dq\alpha_j|^2}} &=&  \frac{\lambda_1^k\dq\alpha_j}{|\lambda_1|^k\sqrt{\sum_{i=1}^l |\dq\alpha_i|^2}}+{\tilde O_D}\left(\left|\frac{\lambda_{l+1,\st}}{\lambda_{1,\st}}\right|^{k}\right)\\
     \nonumber&=& s^k\dq\gamma_j+{\tilde O_D}\left(\left|\frac{\lambda_{l+1,\st}}{\lambda_{1,\st}}\right|^{k}\right).
\end{eqnarray*} 
Here, $\dq\gamma_j=\frac{\dq \alpha_j}{\sqrt{\sum_{i=1}^l |\dq\alpha_i|^2}}$,  the first equality follows from 
     \begin{eqnarray*}
        &&\frac{\lambda_1^k\dq\alpha_j}{\sqrt{\sum_{j=1}^n |\lambda_j^k\dq\alpha_j|^2}} - \frac{\lambda_1^k\dq\alpha_j}{|\lambda_1|^k\sqrt{\sum_{i=1}^l |\dq\alpha_i|^2}}\\
        &=& \frac{O_D(|\lambda_{l+1,\st}|^k)}{|\lambda_1|^{k}\sqrt{\sum_{i=1}^l |\dq\alpha_i|^2}+O_D(|\lambda_{l+1,\st}|^k)}\\
        &=&{\tilde O_D}\left(\left|\frac{\lambda_{l+1,\st}}{\lambda_{1,\st}}\right|^k\right),
    \end{eqnarray*}
    where $\tilde O_D(\cdot)$ is defined by \eqref{equ:titleO_D}, and the second equality follows from $\frac{\lambda_1^k}{|\lambda_1|^k}=s^k$.  
    Furthermore, if  $j\ge l+1$, we have 
\begin{equation*}
    \frac{|\lambda_j^k\dq\alpha_j|}{\sqrt{\sum_{j=1}^n |\lambda_j^k\dq\alpha_j|^2}}=\frac{|\lambda_j|^k|\dq\alpha_j|}{|\lambda_1|^k\sqrt{\sum_{j=1}^l|\dq\alpha_j|^2} + {O_D}(|\lambda_{l+1,\st}|^k)} = 
        {\tilde O_D}\left(\left|\frac{\lambda_{j,\st}}{\lambda_{1,\st}}\right|^k\right). 
\end{equation*}

Consequently,  \eqref{equ:powermethod_vk} holds true and 
$\vdq{v}^{(k)}s^k$ converges to  $\sum_{j=1}^l\vdq u_j \dq \gamma_j$. Since  
\begin{equation*}
    \vdq Q \sum_{j=1}^l\vdq u_j \dq \gamma_j = \sum_{j=1}^l\lambda_j\vdq u_j \dq \gamma_j = \lambda_1\sum_{j=1}^l\vdq u_j \dq \gamma_j,  \end{equation*}
and 
\begin{equation*}
\left\|\sum_{j=1}^l\vdq u_j \dq \gamma_j\right\|_F^2=\left(\sum_{j=1}^l\vdq u_j \dq \gamma_j\right)^*\left(\sum_{j=1}^l\vdq u_j \dq \gamma_j\right)=\sum_{j=1}^l \dq \gamma_j^* \dq \gamma_j=1_D,   
\end{equation*}
we have that $\sum_{j=1}^l\vdq u_j \dq \gamma_j$ is also an eigenvector corresponding to $\lambda_1$.  
Furthermore, 
\begin{equation*}
     \lambda^{(k)} =(\vdq{v}^{(k)})^*\vdq y^{(k+1)} =\frac{\sum_{j=1}^n  \lambda_j^{2k+1}|\dq\alpha_j|^2}{\sum_{j=1}^n |\lambda_j^{k}\dq\alpha_j|^2}= \lambda_1  \left({1_D} + {\tilde O_D}\left(\left|\frac{\lambda_{l+1,\st}}{\lambda_{1,\st}}\right|^{2k}\right)\right).  
\end{equation*} 
This completes the proof. 
\end{proof}


\subsection{General dominant eigenvalues}
In general, we may assume that $\vdq Q$ satisfies Assumption A.    
Then we  show the convergence of  the standard part of the dominant eigenvalue and eigenvector, respectively. 
 
\begin{Lem}
     Given a dual quaternion Hermitian matrix $\vdq Q\in\dqset^{n\times n}$ and a dual quaternion vector $\vdq v^{(0)}=\sum_{j=1}^n\vdq u_j\dq \alpha_j$ with unit norm, where $\vdq u_1,\dots,\vdq u_n$ are orthonormal eigenvectors of $\vdq Q$. Suppose $\vqu Q_{\st}\neq \vqu O_{n\times n}$, $\sum_{j=1}^l |\qu \alpha_{j,\st}|^2\neq0$, and
      \begin{equation*} 
     \lambda_{1,st}=\dots= \lambda_{l,st} \text{ and }|\lambda_{1,\st}|>|\lambda_{l+1,\st}|\ge\dots\ge|\lambda_{n,\st}|\ge0. 
     \end{equation*} 
     Then for the sequence generated by the power method \eqref{equ:DQM_powermethod}, we have 
\begin{equation}\label{equ:g_powermethod_vk}
    \vqu{v}_{\st}^{(k)} = s^k\sum_{j=1}^l\vqu u_{j,\st} \qu \gamma_{j,\st} \left(1 + O\left(\left|\frac{\lambda_{l+1,\st}}{\lambda_{1,\st}}\right|^k\right)\right), 
\end{equation}
where $s=sgn(\lambda_{1,\st})$, $\qu \gamma_{j,\st}=\frac{\qu \alpha_{j,\st}}{\sqrt{\sum_{i=1}^l |\qu\alpha_{i,\st}|^2}}$, and 
\begin{equation*} 
    \lambda_{\st}^{(k)}=
    \left(\vqu{v}_{\st}^{(k)}\right)^*\vqu y_{\st}^{(k+1)} = \lambda_{1,\st} \left(1 + O\left(\left|\frac{\lambda_{l+1,\st}}{\lambda_{1,\st}}\right|^{2k}\right)\right).  
\end{equation*}
 In other words, 
the sequence  $\vqu{v}_{\st}^{(k)}s^k$  converges to  $\sum_{j=1}^l\vqu u_{j,\st} \qu \gamma_{j,\st}$ at a rate of $O\left(\left|\frac{\lambda_{l+1,\st}}{\lambda_{1,\st}}\right|^{k}\right)$,  and the sequence $\lambda_{\st}^{(k)}$ converges to $\lambda_{1,\st}$ at a rate of $O\left(\left|\frac{\lambda_{l+1,\st}}{\lambda_{1,\st}}\right|^{2k}\right)$. 
\end{Lem}
\begin{proof}
 For the standard part of the coefficient  of $\vdq u_j$ in the sequence $\vdq{v}^{(k)}$ presented in  \eqref{equ:vkyk}, we have   
  \begin{equation*}
        \frac{\lambda_{j,\st}^k\qu\alpha_{j,\st}}{\sqrt{\sum_{j=1}^n |\lambda_{j,\st}^k\qu\alpha_{j,\st}|^2}}=\frac{\lambda_{j,\st}^k\qu\alpha_{j,\st}}{|\lambda_{1,\st}|^k\sqrt{\sum_{j=1}^l|\qu\alpha_{j,\st}|^2} + O(|\lambda_{l+1,\st}|^k)}.   
    \end{equation*} 
    If $j\in\{1,\dots,l\}$, it holds that  
\begin{eqnarray*}
    \frac{\lambda_{j,\st}^k\qu\alpha_{j,\st}}{\sqrt{\sum_{j=1}^n |\lambda_{j,\st}^k\qu\alpha_{j,\st}|^2}} &=&  \frac{\lambda_{1,\st}^k\qu\alpha_{j,\st}}{|\lambda_{1,\st}|^k\sqrt{\sum_{j=1}^l |\qu\alpha_{j,\st}|^2}}+O\left(\left|\frac{\lambda_{l+1,\st}}{\lambda_{1,\st}}\right|^{k}\right)\\
     \nonumber&=& s^k\qu\gamma_{j,\st}+O\left(\left|\frac{\lambda_{l+1,\st}}{\lambda_{1,\st}}\right|^{k}\right).
\end{eqnarray*} 
Here, $\qu\gamma_{j,\st}=\frac{\qu \alpha_{j,\st}}{\sqrt{\sum_{i=1}^l |\qu\alpha_{i,\st}|^2}}$. 
If $j\ge l+1$, we have 
\begin{equation*}
   \frac{\lambda_{j,\st}^k\qu\alpha_{j,\st}}{\sqrt{\sum_{j=1}^n |\lambda_{j,\st}^k\qu\alpha_{j,\st}|^2}} =\frac{\lambda_{j,\st}^k\qu\alpha_{j,\st}}{|\lambda_{1,\st}|^k\sqrt{\sum_{j=1}^l |\qu\alpha_{j,\st}|^2} +O(|\lambda_{l+1,\st}|^k)}=
      O\left(\left|\frac{\lambda_{j,\st}}{\lambda_{1,\st}}\right|^{k}\right).
\end{equation*}

Consequently,  \eqref{equ:g_powermethod_vk} holds true and 
    $\vqu{v}_{\st}^{(k)}s^k$ converges to  $\sum_{j=1}^l\vqu u_{j,\st} \qu \gamma_{j,\st}$. Since  
\begin{equation*}
    \vqu Q_{\st} \sum_{j=1}^l\vqu u_{j,\st} \qu \gamma_{j,\st} = \sum_{j=1}^l\lambda_{j,\st}\vqu u_{j,\st} \qu \gamma_{j,\st} = \lambda_{1,\st}\sum_{j=1}^l\vqu u_{j,\st} \qu \gamma_{j,\st},  \end{equation*}
and 
\begin{equation*}
\left\|\sum_{j=1}^l\vqu u_j \qu \gamma_{j,\st}\right\|_F^2=\left(\sum_{j=1}^l\vqu u_{j,\st} \qu \gamma_{j,\st}\right)^*\left(\sum_{j=1}^l\vqu u_{j,\st} \qu \gamma_{j,\st}\right)=\sum_{j=1}^l \qu \gamma_{j,\st}^* \qu \gamma_{j,\st}= 1,   
\end{equation*}
we have $\sum_{j=1}^l\vqu u_{j,\st} \qu \gamma_{j,\st}$ is also an eigenvector of $\vqu Q_{\st}$ corresponding to $\lambda_{1,\st}$.  
Furthermore, $\lambda_{\st}^{(k)}=(\vqu{v}_{\st}^{(k)})^*\vqu Q_{\st}\vqu{v}_{\st}^{(k)}$ converges to $\lambda_{1,\st}$. 
This completes the proof. 
\end{proof} 

\textbf{Remark.} 
In the general case, the dual parts of $\vdq v^{(k)}$ and   $\vdq y^{(k)}$ may not converge. 
In this case,  Algorithm~\ref{alg:power} only returns the standard parts of the strict dominant eigenvalue and its corresponding eigenvector.
Denote the standard part of the eigenvector as $\qu v_{\st}=\sum_{j=1}^l\vqu u_{j,\st} \qu \gamma_{j,\st}$. By \cite{QL22}, we can compute the dual part of the dominant eigenvalue and its corresponding eigenvector  via 
\begin{equation*}
\left\{\begin{array}{l}
      \lambda_{\I}  =  \vqu v_{\st}^*\vqu Q_{\I}\vqu v_{\st},\\   
   \vqu Q_{\st}\vqu v_{\I} + \vqu Q_{\I}\vqu v_{\st} = \lambda_{\st} \vqu v_{\I}+\lambda_{\I} \vqu v_{\st} \text{ and }sc(\vqu v_{\I}^*\vqu v_{\st}) = 0. 
\end{array}
  \right.
\end{equation*}

\bigskip
\subsection{All {appreciable} eigenvalues of a dual quaternion Hermitian matrix} 



By Theorem 7.1 in \cite{LHQ22}, the Eckart-Young-like theorem holds for dual quaternion matrices.    
If $\vdq Q\in\dqset^{n\times n}$ is a dual quaternion  Hermitian  matrix, then by Theorem 4.1 in \cite{QL22},  $\vdq Q$ can be rewritten as 
\begin{equation}\label{equ:DM_EVD}
    \vdq Q=\vdq U  \mathbf \Sigma \vdq U^*=\sum_{i=1}^n  \lambda_i\vdq u_i\vdq u_i^*, 
\end{equation}
where $\vdq U=[\vdq u_1,\dots,\vdq u_n]\in\dqset^{n\times n}$ is a unitary matrix, $ \mathbf \Sigma=\text{diag}(\lambda_1,\dots,\lambda_n)\in\mathbb D^{n\times n}$ is a diagonal dual matrix, and $\lambda_i, i=1,\dots,n$ are in the descending order. 

Denote $\vdq Q_{k}=\vdq Q-\sum_{i=1}^{k-1}\lambda_i \vdq u_i\vdq u_i^*$. Then $\lambda_k,\vdq u_k$ is the dominant eigenpair of $\vdq Q_{k}$. If $\lambda_{k,\st}\neq0$, 
$\lambda_k,\vdq u_k$ can be computed by implementing the power method on $\vdq Q_k$. 
By repeating this process from $k=1$ to $n$, we get all  {appreciable} eigenvalues and their corresponding eigenvectors. 
 

The process is summarized in Algorithm~\ref{alg:power_all}. 
\begin{algorithm}
\caption{Computing all {appreciable}  eigenvalues of a dual quaternion Hermitian matrix}\label{alg:power_all}
\begin{algorithmic}
\Require $\vdq{Q}$,   the dimension $n$, and the tolerance $\gamma$. 
\State $\vdq Q_1  = \vdq Q$.
\For{$k=1,\dots, n$}
 \State Compute $\lambda_k,\vdq u_k$ as the dominant eigenpair of $\vdq Q_k$ by Algorithm \ref{alg:power}. 
 \State Update $\vdq Q_{k+1} = \vdq Q_{k}-\lambda_k\vdq u_k\vdq u_k^*$. 
 \If{{$\|\vqu Q_{k+1,\st}\|_{F}\le \gamma$}}
 \State Stop.
 \EndIf 
 \EndFor
\State \textbf{Output:} Eigenvalues $\{\lambda_i\}_{i=1}^k$ and eigenvectors $\{\vdq u_i\}_{i=1}^k$. 
\end{algorithmic}
\end{algorithm}

\begin{Lem}\label{lem:alleigenvalues}
    Given a dual quaternion Hermitian  matrix $\vdq Q\in\dqset^{n\times n}$. Suppose that all {appreciable}  eigenvalues $\lambda_i$  either equals each other or have distinct standard parts. Then  Algorithm \ref{alg:power_all} can return all  {appreciable}  eigenvalues   and their corresponding eigenvectors. 
\end{Lem}
\begin{proof}
    This lemma follows directly from \eqref{equ:DM_EVD} and Theorem \ref{thm:power_multiplicity_l}. 
\end{proof}


{
\textbf{Remark:} From Lemma~\ref{lem:alleigenvalues} and  Theorem~6.1 in \cite{QL22}, we can also compute the singular value decomposition of a general dual quaternion matrix. 
Specifically, given a dual quaternion matrix  $\vdq A\in\dqset^{m\times n}$, Qi and Luo \cite{QL22} showed that there exist  dual quaternion unitary matrices $\vdq V\in\dqset^{m\times m}$ and $\vdq U\in\dqset^{n\times n}$  such that 
\begin{equation}\label{equ:SVD}
    \vdq A=\vdq V\vdq\Sigma\vdq U^*,
\end{equation}
where $\vdq\Sigma$ is a diagonal matrix. 
Let $\vdq B=\vdq A^*\vdq A$ and $\vdq C=\vdq A\vdq A^*$. Then   
\begin{equation*}
    \vdq B = \vdq U\vdq \Sigma^2\vdq U^* \text{ and } \vdq C = \vdq V\vdq \Sigma^2\vdq V^* 
\end{equation*} 
are both  dual quaternion Hermitian matrices. 
Hence, the singular value decomposition of $\vdq A$ can be obtained by the eigenvalue decompositions of $\vdq B$ and $\vdq C$. 
}
 \subsection{Relationship with the best rank-one approximation}

\begin{Lem}\label{cor:rank1}
Given an {appreciable}  dual quaternion Hermitian  matrix $\vdq Q\in\dqset^{n\times n}$. The dominant eigenvalue and eigenvector formulate the best rank-one approximation of $\vdq Q$ under the {square of the}  $F$-norm, i.e., 
\begin{equation}\label{equ:best_rank1}
    \lambda_1,\vdq u_1 \in\mathop{\arg\min}\limits_{\lambda\in\mathbb D, \vdq u\in\hat{\mathbb Q}_2^n} \quad \|\vdq Q - \lambda \vdq u\vdq u^*\|_F^2. 
\end{equation}
Furthermore, if $\vdq Q$ is rank-one, then $\lambda_1$ and $\vdq u_1$ are also the optimal solutions under the square of the  $F^*$-norm, 
\begin{equation}\label{equ:best_rank_F*}
    \lambda_1,\vdq u_1 \in\mathop{\arg\min}\limits_{\lambda\in\mathbb D, \vdq u\in\hat{\mathbb Q}_2^n} \quad \|\vdq Q - \lambda \vdq u\vdq u^*\|_{F^*}^2. 
\end{equation}
\end{Lem}
\begin{proof}
    Since 
    \begin{eqnarray*}
        &&\|\vdq Q - \lambda \vdq u\vdq u^*\|_F^2=\|\vqu Q_{\st} - \lambda_{\st} \vqu u_{\st}\vqu u_{\st}^*\|_F^2 \\
        &+&2sc\left(tr\left((\vqu Q_{\st} - \lambda_{\st} \vqu u_{\st}\vqu u_{\st}^*)^*(\vqu Q_{\I} - \lambda_{\I} \vqu u_{\st}\vqu u_{\st}^*-\lambda_{\st} \vqu u_{\I}\vqu u_{\st}^*-\lambda_{\st} \vqu u_{\st}\vqu u_{\I}^*)\right)\right)\epsilon 
    \end{eqnarray*}
    and $\vqu Q_{\st}$ is a quaternion Hermitian matrix, we have 
    \begin{equation*}
       \lambda_{1,\st},\vqu u_{1,\st} \in \mathop{\arg\min}\limits_{\lambda\in\mathbb R,\vqu u\in\mathbb Q_2^{n\times 1}}\|\vqu Q_{\st} - \lambda  \vqu u \vqu u^*\|_F^2.
    \end{equation*}
    Further, by $\vqu Q_{\st}\vqu u_{1,\st}= \lambda_{1,\st} \vqu u_{1,\st}$,  $\vqu u_{1,\st}^*\vqu u_{1,\st}=\qu 1$, and $sc(\vqu u_{1,\st}^*\vqu u_{\I})= 0$, we have 
    \begin{eqnarray*}
        &&sc\left(tr\left((\vqu Q_{\st} - \lambda_{1,\st} \vqu u_{1,\st}\vqu u_{1,\st}^*)^*(\vqu Q_{\I} - \lambda_{\I} \vqu u_{1,\st}\vqu u_{1,\st}^*-\lambda_{1,\st} \vqu u_{\I}\vqu u_{1,\st}^*-\lambda_{1,\st} \vqu u_{1,\st}\vqu u_{\I}^*)\right)\right)\\
        &=& sc(tr(\vqu Q_{\st}\vqu Q_{\I})) -  \lambda_{1,\st}\lambda_{\I}-  \lambda_{1,\st}^2 sc(\vqu u_{1,\st}^*\vqu u_{\I} + \vqu u_{\I}^*\vqu u_{1,\st})\\
        && - \lambda_{1,\st} \vqu u_{1,\st}^*\vqu Q_{\I}\vqu u_{1,\st}  +  \lambda_{1,\st}\lambda_{\I}+\lambda_{1,\st}^2 sc(\vqu u_{1,\st}^*\vqu u_{\I} + \vqu u_{\I}^*\vqu u_{1,\st})\\
        & = &sc(tr(\vqu Q_{\st}\vqu Q_{\I})) - \lambda_{1,\st} \vqu u_{1,\st}^*\vqu Q_{\I}\vqu u_{1,\st}  
    \end{eqnarray*}
    for any $\vqu u_{\I}$ satisfying $sc(\vqu u_{1,\st}^*\vqu u_{\I})= 0$. 
    Hence,   $\lambda_1$ and $\vdq u_1$ is an  optimal solution of \eqref{equ:best_rank_F*}.

   In addition, if $\vdq Q$ is a rank-one matrix, then $\vdq Q=\lambda_1\vdq u_1\vdq u_1^*$ and 
    \begin{eqnarray*}
         &&\|\vdq Q - \lambda_1 \vdq u_1\vdq u_1^*\|_{F*}^2  =\|\vqu Q_{\st} - \lambda_{1,\st} \vqu u_{1,\st}\vqu u_{1,\st}^*\|_F^2=0_D.
    \end{eqnarray*}
    Hence   $\lambda_1,\vdq u_1$ is   the optimal solution of \eqref{equ:best_rank_F*} as it attains the global optimal value. This completes the proof. 
\end{proof}

\bigskip 
\section{Numerical Experiments for Computing Eigenvalues}\label{sec:numerical}
In this section, we first show the numerical experiments for computing all {appreciable} eigenvalues of the Laplacian matrices of   circles and   random graphs.

\subsection{Eigenvalues of Laplacian matrices of circles}
In the multi-agent formation control, the Laplacian matrix of the mutual visibility graph plays a key role \cite{QWL22}.  
Given a vector $\vdq q\in\hat{\mathbb U}^{n\times 1}$ and a graph $G=(V,E)$ with $n$ vertices and $m$ edges, the Laplacian matrix of $G$ is defined as follows, 
\begin{equation*}
    \vdq L =  \mathbf D-\vdq A, 
\end{equation*}
where $ \mathbf D$ is a diagonal real matrix where the $i$-th diagonal element  is equal to the degree of  the $i$-th vertex, and $\vdq A=(\dq a_{ij})$, 
\begin{equation*}
    \dq a_{ij} = \left\{
\begin{array}{cl}
   \dq q_i^*\dq q_j,  &  \text{if } (i,j)\in E,\\
    0, & \text{otherwise}. 
\end{array}
    \right.
\end{equation*} 
In the multi-agent formation control, $\dq q_i\in\hat{\mathbb U}$ is the   configurations of the $i$-th  rigid body,  $\dq a_{ij}\in\hat{\mathbb U}$ is the relative configuration of the rigid bodies $i$ and $j$, and $\vdq A$ is  
the relative configuration adjacency matrix. 

\begin{figure}
    \centering
    \includegraphics[width=0.4\textwidth]{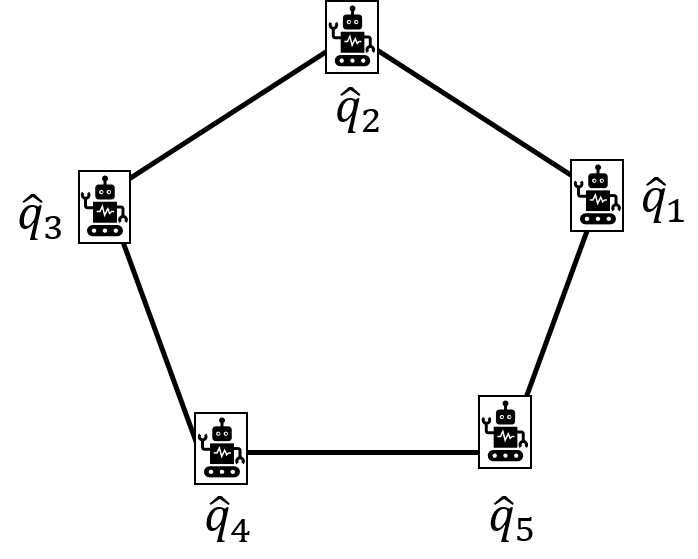}
    \caption{A five-point  circle}
    \label{fig:5pcircle}
\end{figure}
Consider a five-point circle as shown in Fig.~\ref{fig:5pcircle}. Suppose $\vdq q$ is a  random  unit dual quaternion vector   as follows, 
{\small \begin{equation*}
   \vdq q = \left[\begin{array}{rrrr}
   -0.4568 &  -0.0602 &   -0.0555  &  -0.8858 \\
    0.6545 &   -0.5512 &    0.0023 &   -0.5175 \\
    0.3086  &   0.7100 &   -0.5253 &   -0.3533 \\
   -0.5730  &   0.1223 &   -0.6611  &  -0.4688 \\
   -0.5851 &    0.0650  &  -0.1431  &  -0.7956  \\ 
    \end{array}
    \right]
    +\left[\begin{array}{rrrr}
    0.4701  &  -0.6467   &  0.7286  &  -0.2441 \\
    0.1108   & -0.4486  &   0.6871  &   0.6210 \\
   -0.6448   & -0.1852   & -0.2048  &  -0.6311 \\
   -0.8584  &  -0.2427   &  1.2512  &  -0.7785 \\
   -0.2806  &   0.4410   &  0.1730  &   0.2113 \\
    \end{array}
    \right] \epsilon. 
\end{equation*}} 

\begin{figure}
\centering 		
\subfloat[The first eigenpair for $n=5$.]{ 			
\label{fig:power_convergence_1}			
\includegraphics[scale=0.25]{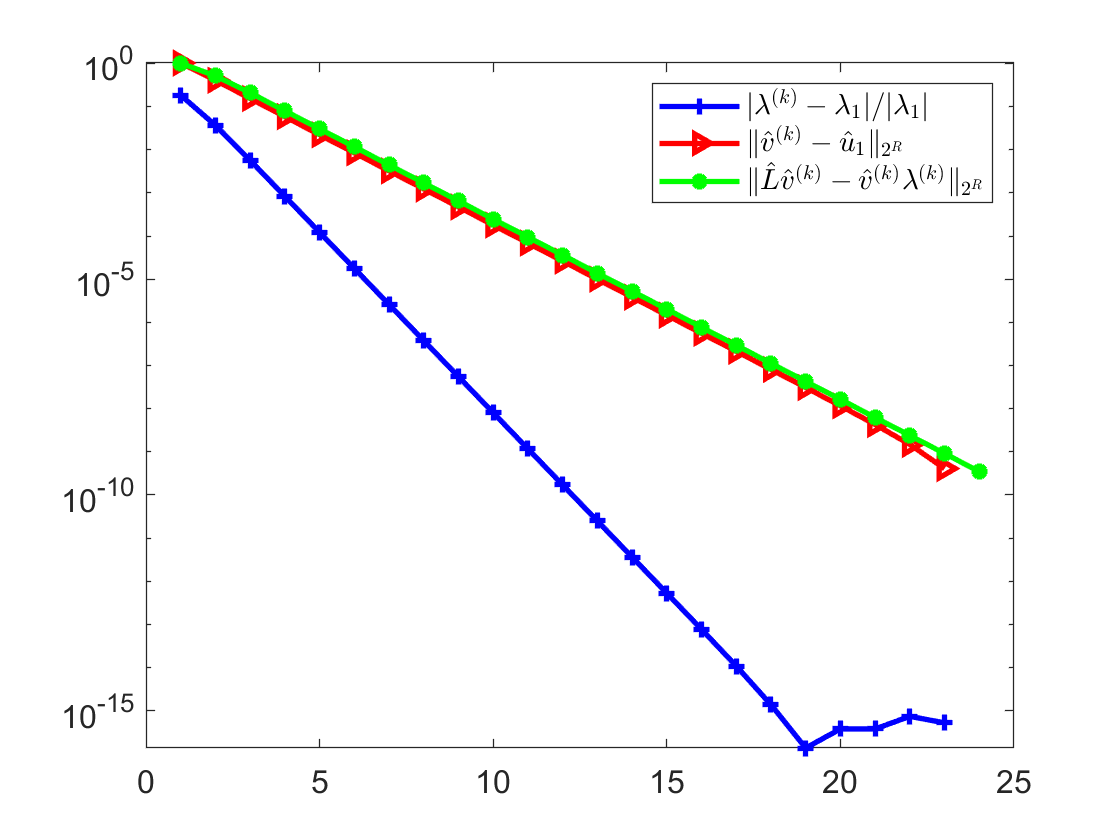}} 		
\hspace{0.1in} 
\subfloat[The second eigenpair for $n=5$.]{ 			
\label{fig:power_convergence_2}  			
\includegraphics[scale=0.25]{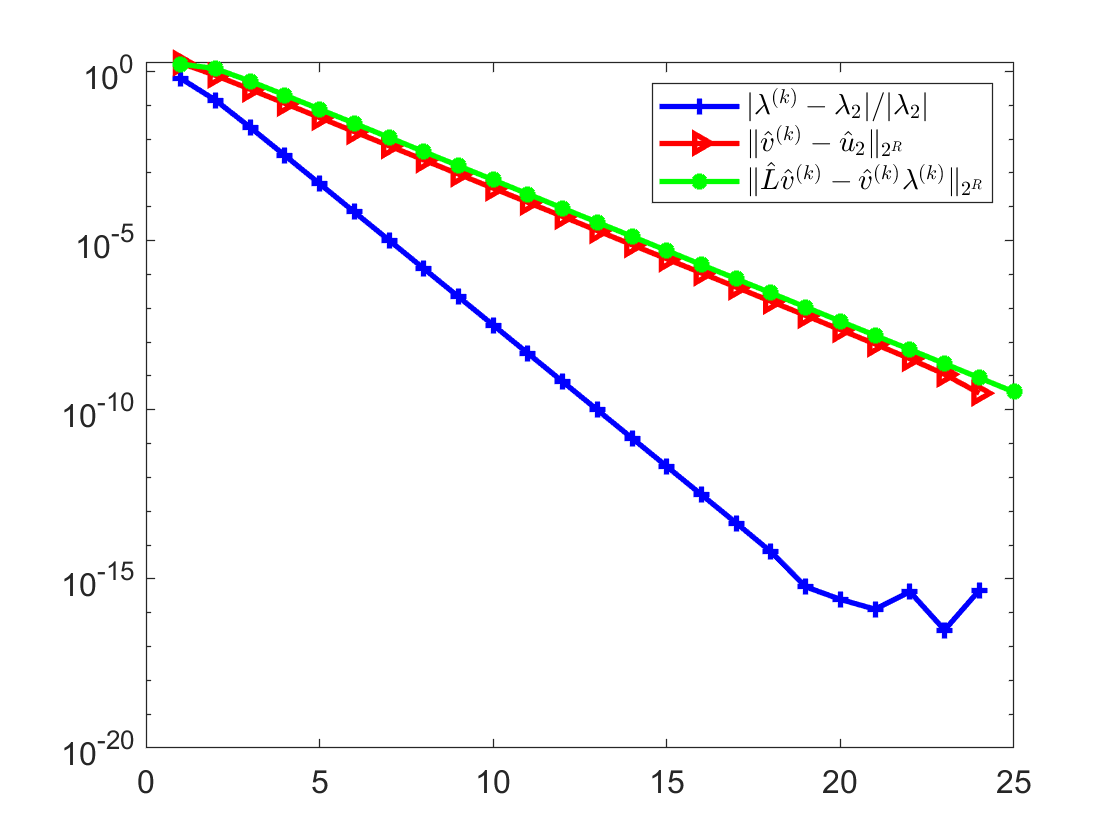}} 
\caption{Convergence rate of power method for the Laplacian matrix of a five-point circle.} 		
\label{fig:power_convergence}  	
\end{figure}

We now compute all appreciable eigenvalues of $\vdq L$ by Algorithm~\ref{alg:power_all}. After 23, 25, 1, and 1 iterations respectively, we get four  eigenvalues 3.618, 3.618, 1.382, and 1.382.  Then Algorithm~\ref{alg:power_all} is stopped since the stopping criterion $\|\vqu L_{5,\st}\|_F\le 10^{-6}\times \|\vqu L_{\st}\|_F$ is satisfied. Here, $\vdq L_5=\vdq L-\sum_{k=1}^4\lambda_k\vdq u_k\vdq u_k^*$. Further,  $\vqu L_{5,\I}$ is also close to a zero matrix. Hence, all five eigenvalues of $\vdq L$ are:   
\begin{equation*}
     [3.6180+0\epsilon,\  3.6180+0\epsilon, \  1.3820+0\epsilon,\   1.3820+0\epsilon,\   0+0\epsilon].
\end{equation*}
 The iterates of $\lambda^{(k)}$  and $\vdq v^{(k)}$  for computing the first and second eigenpair are shown in Figure~\ref{fig:power_convergence}. 
From this figure, we see that both the eigenvalue and the eigenvector converge  linearly and the  eigenvalue converges  much faster than the eigenvector. 
Comparing the first and the second eigenpair, we see the convergence rate  is the same. These conclusions are corresponding to our theory in Theorem~\ref{thm:power_multiplicity_l} that  $\vdq v^{(k)}$ converges to a strict dominant eigenvector   at a rate of $\tilde O_D((\frac{1.382}{3.618})^k)$ and $\lambda^{(k)}$ converges to a strict dominant eigenvalue   at a rate of $\tilde O_D((\frac{1.382}{3.618})^{2k})$. 


\begin{figure}
\centering 		
\subfloat[The first eigenpair for $n=6$.]{ 			
\label{fig:power_convergence_1_n6}			
\includegraphics[scale=0.25]{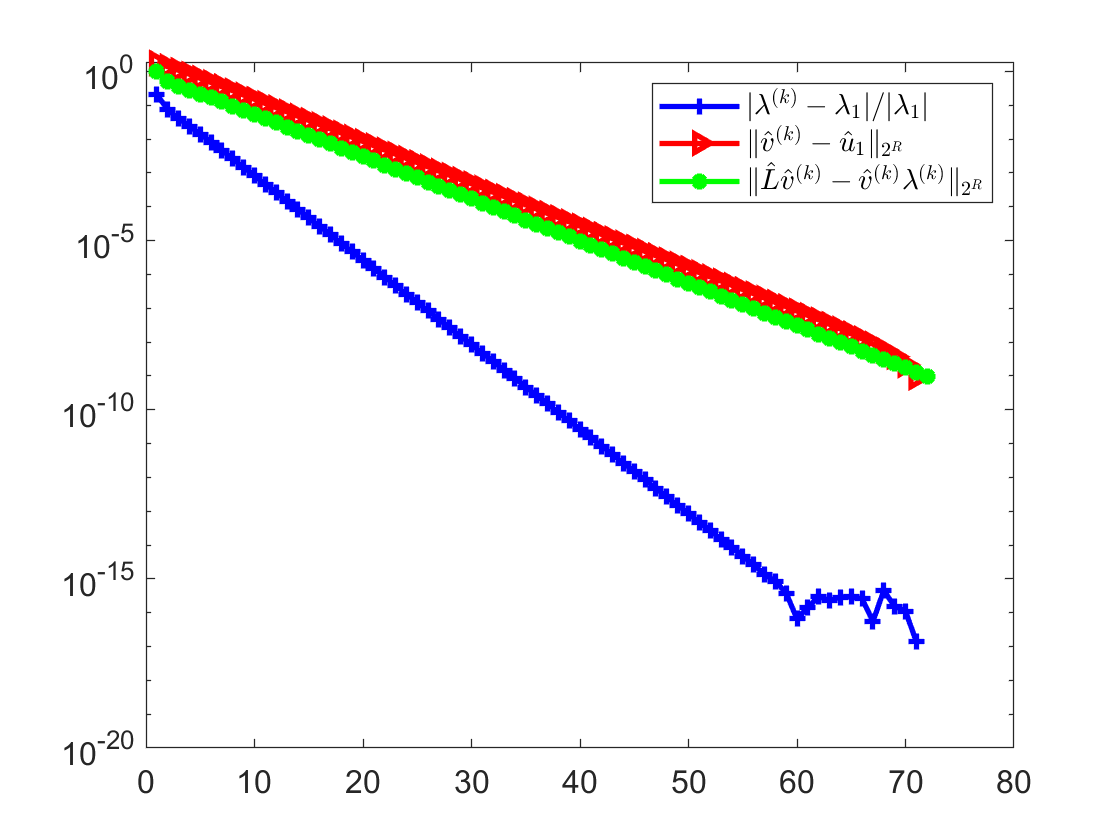}} 		
\hspace{0.1in} 
\subfloat[The second eigenpair for $n=6$.]{ 			
\label{fig:power_convergence_2_n6}  			
\includegraphics[scale=0.25]{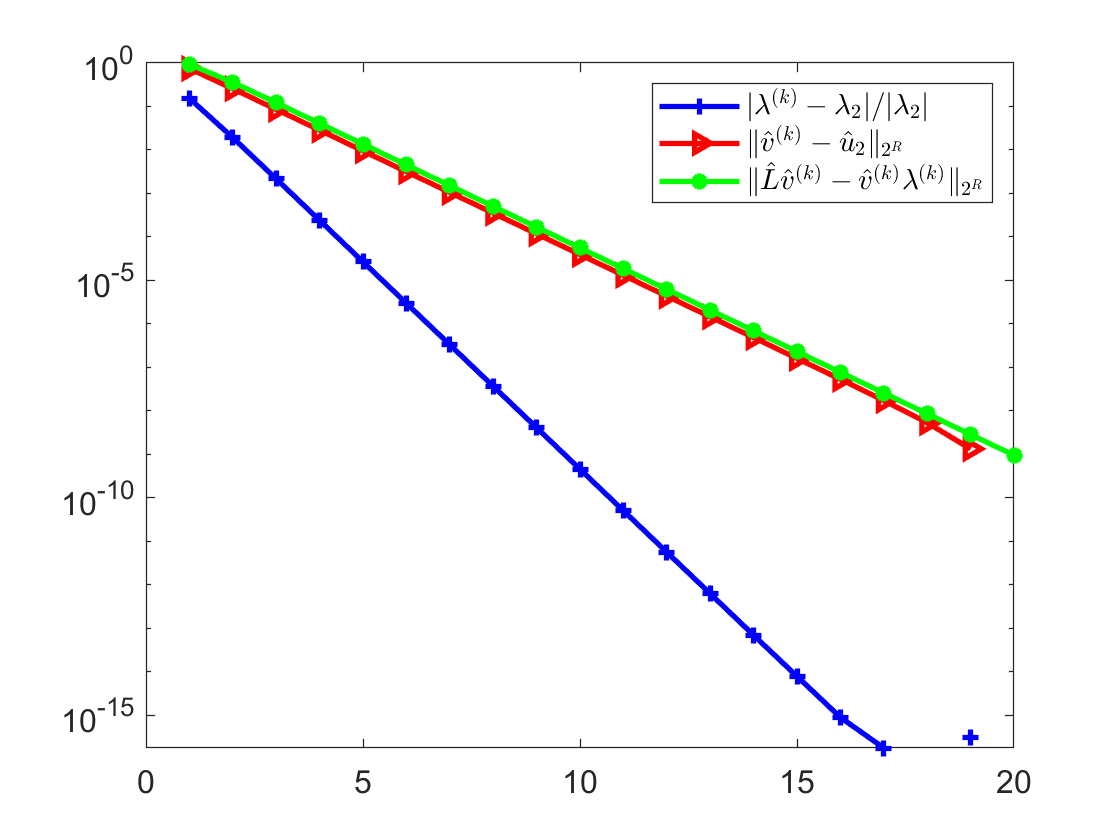}} 
\caption{Convergence rate of power method for the Laplacian matrix of a six-point circle.} 		
\label{fig:power_convergence_n6}  	
\end{figure}

We further consider the Laplacian matrix of  circles with   3 to 10 points and   list their eigenvalues  and the number of iterations of  the power method  in  Table~\ref{tab:eigs_5pcircle_Laplacian}.   
From this table, we see that all {appreciable} eigenvalues are real numbers and the smallest eigenvalues are all zeros. Further, the second smallest eigenvalue is monotonically decreasing with the dimension $n$. 
{This is due to the fact that the second smallest eigenvalue   is corresponding to the algebra connectivity of the graph.} 
For each $n$, we see the eigenvalues with the same values consume almost the same number of iterations. This is corresponding to our theory that the convergence rate is $\tilde O_D\left(\left|\frac{\lambda_{l+1,\st}}{\lambda_{1,\st}}\right|^{2k}\right)$. 
Further, we see that the number of iterations is less   when $\left|\frac{\lambda_{l+1,\st}}{\lambda_{1,\st}}\right|$ is small. 
We  also present the convergence rate of $\lambda^{(k)}$  and $\vdq v^{(k)}$  of the first and second eigenpairs for the six-point circle  in Figure~\ref{fig:power_convergence_n6}. 
We see that the sequences  converge linearly. The power method for computing the first and the second  eigenvalues converge in 76 and 20 iterations, respectively. This is corresponding to our theory that  the sequences for computing the first and the second eigenvalues converge at a rate of $\tilde O_D\left(\left(\frac{3}{4}\right)^{2k}\right)$ and $\tilde O_D\left(\left(\frac{1}{3}\right)^{2k}\right)$, respectively.

\begin{table}
    \centering
     \caption{The eigenvalues (top row)   for the Laplacian matrices of the circle and and the number of iterations (bottom row) of the power method.} 
    \begin{tabular}{|c|cccccccccc|}
    \hline
    $n$ & $[\lambda_1,$&$\lambda_2,$&$\lambda_3,$&$\lambda_4,$&$\lambda_5,$&$\lambda_6,$&$\lambda_7,$&$\lambda_8,$&$\lambda_9,$&$\lambda_{10}]$\\
    \hline
        $3$ &  [3.0000,  &   3.0000,  &   0]&&&&&&&\\
        & [1,& 1,& 0]&&&&&&&\\ \hline
         $4$ & [4.0000,  &  2.0000,  &  2.0000,   &      0] &&&&&& \\
         & [33, &1,& 1,& 0]  &&&&&&\\ \hline
         $5$ & [3.6180, &   3.6180,   & 1.3820, &   1.3820,   &      0]&&&&&  \\
          & [24,& 25, &1, & 1, & 0] &&&&&\\\hline
         $6$ & [4.0000,  &  3.0000,   & 3.0000, &   1.0000,  &  1.0000,    &     0] &&&& \\
          & [76, &20, &20,  &1, & 1,  &0]&&&& \\\hline
         7 & [3.8019,  & 3.8019, &   2.4450,  &  2.4450,  &  0.7530,   & 0.7530,       &  0]&&& \\
          & [50, & 48, & 20, & 20, & 1, & 1, & 0] &&&\\\hline
         8 & [4.0000,   & 3.4142,  &  3.4142,   & 2.0000,  &  2.0000,  &   0.5858,  &  0.5858,     &    0]&& \\
          & [133,  &  41,  &  41,  & 18, &   20, &1, &1,& 0] && \\\hline
         9 & [3.8794,  &  3.8794,  &  3.0000,  &  3.0000,   & 1.6527, &   1.6527,&    0.4679,  &  0.4679,     &    0]& \\
          & [ 81,  &  85,  &  37, &  39,  &  18,  &  19, & 1,& 1,& 0] & \\\hline
         10 & [4.0000, &   3.6180,  &  3.6180,  &  2.6180 ,&   2.6180,  &  1.3820,  &  1.3820, &   0.3820,   & 0.3820,    &     0] \\
          & [203,   & 67,  &  70,   & 35, &   36,  &  18,  &  18,& 1, &1,& 0] \\
         \hline
    \end{tabular} 
    \label{tab:eigs_5pcircle_Laplacian}
\end{table}

In our numerical experiments, we find  that the eigenvalues are the same for any $\vdq q$. Hence, we make the following conjecture. 

\textbf{Conjecture.} The eigenvalues of the Laplacian matrix of a circle are all nonnegative real values, and they are independent of the choice of $\vdq q$. 

\subsection{Eigenvalues of Laplacian matrices of random graphs}
We continue to consider the Laplacian matrix of a random graph $G=(V,E)$. Assume that $G$ is a sparse undirected graph  and $E$ is symmetric. The sparsity of $G$ is $s=m/n^2$, where $m=|E|$ is the number of edges in $E$. In practice, we randomly generate $\lceil\frac m2\rceil$ edges and let $E=\{(i,j),(j,i):\text{ if }(i,j) \text{ is sampled}\}$.  
We show the results  with $n=10$ and $n=100$ in Table~\ref{tab:Laplacian_random}.  
All results are repeated ten times with different choices of $\vdq q$, different $E$, and different initial values in the power method and the average results are reported. Denote the error in eigenvalues and the residue of the result matrix as 
\begin{equation*}
    e_{\lambda} = \frac{1}{n_0}\sum_{i=1}^{n_0}\|\vdq L \vdq u_i-\lambda_i \vdq u_i\|_{2^R}  \text{ and } e_L=\frac1{\|\vdq L\|_{F^R}}\|\vdq L  -\sum_{i=1}^{n_0}\lambda_i \vdq u_i\vdq u_i^*\|_{F^R}, 
\end{equation*}
where $\lambda_i$ and $\vdq u_i$ are the limit points of $\lambda^{(k)}$ and $\vdq v^{(k)}$ for computing $\vdq L_i$ in Algorithm \ref{alg:power} and $n_0$ is the total number of eigenvalues obtained in Algorithm~\ref{alg:power_all}. 
In Table~\ref{tab:Laplacian_random}, we denote 
`$n_{\text{iter}}$' as the average   number of iterations  and `time (s)'  as the average CPU time  in seconds for computing one  eigenvalue.

From Table~\ref{tab:Laplacian_random}, we see that when $n=10$, $e_{\lambda}$ is less than $10^{-7}$ and $e_L$ is less than $10^{-14}$. 
Similarly, when $n=100$,  the error in eigenvalues  and the residue of the result matrices is less than $1.1\times 10^{-3}$ and $10^{-7}$, respectively.  
All results are obtained in less than 1.1 seconds. This shows the efficiency of the power method in computing all {appreciable} eigenvalues. 

{\small \begin{table}
    \centering 
    \caption{Numerical results of power method for computing all {appreciable} eigenvalues of Laplacian matrices  of random graphs with different sparsity values.}  
    \begin{tabular}{|c|cccc||c|cccc|}
    \hline
    \multicolumn{10}{|c|}{$n=10$} \\
    \hline 
    $s$ & $e_{\lambda}$ & $e_L$ & $n_{\text{iter}}$ & time (s) & 
     $s$ & $e_{\lambda}$ & $e_L$ & $n_{\text{iter}}$ & time (s)    \\ 
     \hline

         10\% &   2.27e$-$10   &   1.93e$-$16   &          13.72   &     0.0088 &  
          20\% &   4.25e$-$10   &   4.03e$-$16   &        45.25  &        0.0288\\ 
          30\% &   1.31e$-$9  &    4.24e$-$16   &        88.98  &       0.0560 &  
          40\% &   9.31e$-$10  &    2.40e$-$16  &        100.2  &       0.0631 \\
          50\% &   1.72e$-$9  &     2.91e$-$16  &        120.63  &       0.0763 &  
          60\% &   4.43e$-$8  &    4.98e$-$15   &       159.12   &      0.1005 \\  
    \hline    \multicolumn{10}{|c|}{$n=100$} \\ \hline 
    
          5\% &     1.42e$-$4 &  6.27e$-$8    &       707.85  &      0.7880 & 
         8\% &     2.26e$-$4   &   2.38e$-$8    &          797.20  &      0.8857\\
          10\% &     2.99e$-$4 &     2.80e$-$8  &       840.98   &       0.9520&  
          15\% &     7.92e$-$4  & 2.18e$-$8   &      892.15   &     1.078 \\
          18\% &     8.86e$-$4 &  1.45e$-$8   &      913.18   &     1.075 &
          20\% &    1.10e$-$3 &  3.91e$-$8    &   921.19   &     1.071 \\
    \hline 
    \end{tabular}
    \label{tab:Laplacian_random}
\end{table}
}

\bigskip

  \section{The Dominant Eigenvalue Method   for SLAM}\label{sec:SLAM}
    {In this section, we first reformulate the SLAM problem as a rank-one dual quaternion matrix recovery problem. Then we present a two-block coordinate approach for solving this model. The first block has a closed-form solution  and the second block is a rank-one  dual quaternion matrix decomposition problem whose solution can be obtained by the dominant eigenpair. At last, we present several numerical experiments to show the efficiency of our proposed method.}  
    
      Given a directed graph $G = (V, E)$ with $n$ vertices and $m$ edges,  where each vertex $i \in V$ corresponds to a robot pose $\dq q_i \in {\hat {\mathbb U}}$ for $i = 1, \cdots, n$, and each directed edge (arc) $(i, j) \in E$ corresponds to a relative measurement $\dq q_{ij}\in  \hat {\mathbb U}$.     The aim is to find the best $\dq q_i$ for $i = 1, \cdots, n$, to satisfy \cite{CTDD15}
  \begin{equation}\label{equ:SLAM}
   \dq q_{ij} = (\dq q_i)^* \dq q_j.
  \end{equation}
  for $(i, j) \in E$.

\subsection{Reformulation} 
  Denote $\vdq x=[\dq q_1^*,\dots,\dq q_n^*]\in \hat {\mathbb Q}^{n\times 1}$ as a dual quaternion vector, $\vdq Q=(\dq q_{ij})\in \hat {\mathbb Q}^{n\times n}$ as a dual quaternion matrix.  Then $\vdq x\in\hat{\mathbb U}^{n\times 1}$  is a unit vector and $\vdq Q\in\hat{\mathbb U}^{n\times n}$ is a Hermitian matrix.    Equation \eqref{equ:SLAM} is equivalent to  
    \begin{equation}\label{equ:SLAM_M}
   (\vdq Q)_{E} = (\vdq x  \vdq x^*)_E.
  \end{equation}
  Here, $(\vdq Q)_{E}$ denotes the set $\{\dq q_{ij}: (i,j)\in E\}$.

 The least square model aims to find $\vdq x$   as follows, 
   \begin{equation}\label{equ:SLAM_LS}
   \min_{\vdq x \in \hat {\mathbb U}^{n\times 1}}\quad \|P_E(\vdq x  \vdq x^*-\vdq Q)\|_F^2.
  \end{equation}
    The solution of \eqref{equ:SLAM_LS} is not unique because $\vdq x  \vdq x^*=(\vdq x \dq q)  (\vdq x\dq q)^*$ and $\vdq x \dq q  \in \hat {\mathbb U}^{n\times 1}$ for any $\dq q\in\hat{\mathbb U}$. 
    
        

Instead of solving \eqref{equ:SLAM_LS} directly, 
we introduce the auxiliary variables $\vdq X$ to solve the following problem 
  \begin{eqnarray}\label{equ:SLAM_aux}
     &\min\limits_{\vdq X\in \hat{\mathbb Q}^{n\times n}}&\quad  \frac12\|P_E(\vdq X-\vdq Q)\|_F^2 \\
  \nonumber & \text{s.t.} & \vdq X\in C_1 \cap C_2. 
  \end{eqnarray} 
Here, 
\begin{eqnarray}
    C_1 &=& \{\vdq X\in \hat{\mathbb Q}^{n\times n}:\ \vdq X = \vdq X^*, \dq x_{ij}\in\hat{\mathbb U}, \text{ and } \dq x_{ii}=\hat 1\},\\ 
C_2&=&  \{\vdq X\in \hat{\mathbb Q}^{n\times n}:\  \rank(\vdq X)=1\}.
\end{eqnarray} 
 In the following, we   show the equivalence between problems \eqref{equ:SLAM_LS} and \eqref{equ:SLAM_aux}.  
 
  \begin{Thm}\label{thm:equivalence} 
      The optimization problems \eqref{equ:SLAM_LS} and \eqref{equ:SLAM_aux} are equivalent in the following senses: 
      \begin{itemize}
          \item[(i)] If $\vdq x$ is an optimal solution of \eqref{equ:SLAM_LS}, then $\vdq X =   \vdq x \vdq x^*$  is an optimal solution of \eqref{equ:SLAM_aux}. 
          Furthermore, the optimal value of problems \eqref{equ:SLAM_LS} and \eqref{equ:SLAM_aux} are the same. 

          \item[(ii)] Let  $\vdq X$ be the optimal solution for \eqref{equ:SLAM_aux}. Then its rank-one decomposition $\vdq X=\lambda \vdq u\vdq u^*$ satisfies  $\lambda>0_D$. 
          Furthermore,    $\vdq x = \sqrt{\lambda}\vdq u$ is an optimal solution of \eqref{equ:SLAM_LS}, and  
           the optimal value of problems \eqref{equ:SLAM_LS} and \eqref{equ:SLAM_aux} are the same. 
      \end{itemize}
  \end{Thm}
  \begin{proof}
  
  We first show that there is a one-to-one correspondence between $\vdq x\in\hat{\mathbb U}^{n\times 1}$  and $\vdq X\in C_1\cap C_2$.  
  On one hand, suppose that $\vdq x\in\hat{\mathbb U}^{n\times 1}$. Then   it follows from direct derivations that 
  \begin{equation}\label{equ:Xxx}
  \vdq X=\vdq x\vdq x^*  \in C_1\cap C_2.
  \end{equation} 
  On the other hand, suppose that $\vdq X\in C_1\cap C_2$.  Then  it follows from \eqref{equ:DM_EVD} that  there exists  $\lambda \in \mathbb D$ and $\vdq u^*\vdq u=\hat 1$ such that $\vdq X =\lambda\vdq u\vdq u^*$. 
      By $x_{ii} = \lambda\dq u_i \dq u_i^* =\hat 1,$  we have 
      $\lambda_{\st}\qu  u_{i,\st}\qu u_{i,\st}^* =\tilde 1.$ 
     Hence, $\lambda_{\st}>0$ and $\lambda> 0_D$. By \cite{QLY22},  we have $\sqrt{\lambda}=\sqrt{\lambda_{\st}} + \frac{\lambda_{\I}}{2\sqrt{\lambda_{\st}}}\epsilon$. 
      Let $\vdq x=\sqrt{\lambda}\vdq u$. Then $\dq x_i^*\dq x_i=\dq 1$. Hence, 
      $\vdq x\in\hat{\mathbb U}^{n\times 1}$. 
     
      $(i)$  From \eqref{equ:Xxx} we have  $\vdq X =   \vdq x \vdq x^*$  is a feasible solution of  \eqref{equ:SLAM_aux}. 
      Suppose that $\vdq X$ is  not  an optimal solution. Denote the optimal solution as $\vdq X'$.  Then $\|P_E(\vdq X'-\vdq Q)\|_F^2<\|P_E(\vdq X-\vdq Q)\|_F^2$. It follows from the above discussion, there also exists $\vdq x'\in \hat{\mathbb U}^{n\times 1}$ such that 
     $\|P_E(\vdq x'(\vdq x')^*-\vdq Q)\|_F^2<\|P_E(\vdq x  \vdq x ^*-\vdq Q)\|_F^2$. This contradicts   the optimality of $\vdq x$.

      $(ii)$  Denote  $\vdq x = \sqrt{\lambda}\vdq u$. Then it is a feasible  solution of \eqref{equ:SLAM_LS}. 
      The optimality of $\vdq x$ can be obtained by the contradiction method similarly.  
      
      This completes the proof. 
  \end{proof}




\subsection{A block coordinate descent method}
We introduce the auxiliary variables $\vdq X_1$ and $\vdq X_2$ to solve the following problem 
  \begin{eqnarray}\label{equ:SLAM_aux2b}
     &\min\limits_{\{\vdq X_i\}_{i=1}^2\in \hat{\mathbb Q}^{n\times n}}&\quad  \frac12\|P_E(\vdq X_1-\vdq Q)\|_F^2 \\
  \nonumber & \text{s.t.} & \vdq X_1\in C_1, \ \vdq X_2\in C_2,\ \vdq X_1=\vdq X_2. 
  \end{eqnarray} 
The quadratic penalty approach is applied to reformulate \eqref{equ:SLAM_aux2b} as  
    \begin{eqnarray}
      &\min\limits_{\{\vdq X_i\}_{i=1}^2\in \hat{\mathbb Q}^{n\times n}}&\quad  \frac12\|P_E(\vdq X_1-\vdq Q)\|_F^2 +\frac{\rho}{2}\|\vdq X_1-\vdq X_2\|_F^2,\label{equ:SLAM_aux2}\\
  \nonumber & \text{s.t.} & \vdq X_1\in C_1,  \vdq X_2\in C_2. 
  \end{eqnarray}   
Then   we compute $\vdq X_1$ and $\vdq X_2$    alternatively by the block coordinate descent method \cite{XY13}.  
  
  Specifically, at the $k$-th iteration, given $\rho^{(k-1)}$, $\vdq X_1^{(k-1)}$ and $\vdq X_2^{(k-1)}$,  we compute $\vdq X_1^{(k)}$ as follows,
  \begin{eqnarray*}
      \vdq X_1^{(k)}= \mathop{\arg\min}_{\vdq X_1 \in C_1} \ \frac12\|P_E(\vdq X_1 -\vdq Q)\|_F^2 + \frac{\rho^{(k-1)}}{2}\|\vdq X_{1} -\vdq X_{2}^{(k-1)}\|_F^2,   
  \end{eqnarray*}
    which has an explicit solution as $\dq x_{1,ii}^{(k)}=1$ and 
    \begin{equation}\label{equ:X1}
      \dq x_{1,ij}^{(k)} =  
     P_{\hat{\mathbb U}}\left(c_{ij}\left(\delta_{E,ij} \dq q_{ij}+\delta_{E,ji}\dq q_{ji}+\rho^{(k-1)}( \dq x_{2,ij}^{(k-1)} + \dq x_{2,ji}^{(k-1)} )\right)\right)   
  \end{equation}
 for all $i\neq j$. 
Here, $c_{ij}=\frac{1}{\delta_{E,ij}+\delta_{E,ji}+2\rho^{(k-1)}}$,  $\delta_{E,ij}$ is equal to one if $(i,j)\in E$ and zero otherwise, and $P_{\hat{\mathbb U}}(\cdot)$ is the projection  onto the set of unit dual quaternion numbers  defined by \eqref{equ:dq_normalize}.

Then we compute  $\vdq X_2^{(k)}$ as follows,
  \begin{equation}\label{equ:X2subprob}
      \vdq X_2^{(k)}= \mathop{\arg\min}_{\vdq X_2\in C_2} \quad  \frac{\rho^{(k-1)}}{2}\|\vdq X_{2}-\vdq X_{1}^{(k)}\|_F^2. 
  \end{equation}
  {By \eqref{equ:best_rank1},}  this problem admits an explicit solution as 
    \begin{equation}\label{equ:X2}
      \vdq X_2^{(k)} = \lambda\vdq u\vdq u^*, 
  \end{equation}
  where $\lambda$ and $\vdq u$ are the dominant eigenvalue and eigenvector of $\vdq X_1^{(k)}$, respectively, which can be obtained by the power method \eqref{equ:DQM_powermethod}.

 We update the penalty parameter by 
 \begin{equation}\label{equ:rho}
     \rho^{(k)} = \rho_1\rho^{(k-1)},
 \end{equation}
 where $\rho_1>1$ is a constant. 

  We summarize the block coordinate descent method for computing the SLAM problem in Algorithm~\ref{alg:SLAM}. 
  
\begin{algorithm}[t]
\caption{A two-block coordinate descent algorithm for  SLAM}\label{alg:SLAM}
\begin{algorithmic}
\Require $(\vdq{Q})_E$,  the initial values  $\rho^{(0)}$, $\vdq X_1^{(0)}$,   $\vdq X_2^{(0)}$, the maximal iteration $k_{\max}$, the parameter $\rho_1$, and the tolerance $\beta$. 
\For{$k=1,\dots, k_{\max}$}
 \State Update $\vdq X_1^{(k)}$ by \eqref{equ:X1}.
 \State Update $\vdq X_2^{(k)}$ by \eqref{equ:X2}.
 \State Update $\rho^{(k)}$ by \eqref{equ:rho}.  
 \If{$\|\vdq X_1^{(k)}-\vdq X_2^{(k)}\|_{F^R}\le \beta$}
 \State Stop.
 \EndIf
 \EndFor
\State \textbf{Output:} $\vdq X_1^{(k)}$   and $\vdq X_2^{(k)}$.
\end{algorithmic}
\end{algorithm}

\subsection{Numerical results for SLAM}

   We now consider solving the SLAM problem by our proposed model \eqref{equ:SLAM_aux}.    
   Suppose there are $m$ directed edges in $E$.   
  Denote   $\vdq Q_0=\vdq x\vdq x^*\in \hat {\mathbb U}^{n\times n}$ and 
  \begin{equation}\label{SLAM:noise}
      \vdq Q = \vdq Q_0+\vdq N,
  \end{equation}
 where $\vdq N$ is the noise matrix. 
  

The parameters of Algorithm~\ref{alg:SLAM} are given as follows. The penalty parameters $\rho^{(0)}=0.01$ and $\rho_1=1.1$. The maximal iteration is set to be 1000, and the parameter in the stopping criterion is $\beta=10^{-5}$.  
 All experiments are repeated ten times with different choices of $\vdq q$, different noises, and different initial values in Algorithm~\ref{alg:SLAM}.  We only report the average performance. 

Firstly, consider the five-point circle SLAM problem with noisy observation  \eqref{SLAM:noise}.  
The results for different noise levels are given in Table~\ref{tab:SLAM_5pcircle}.
Here, `$l_{\text{noise}}$' is the relative noise level defined by $l_{\text{noise}}=\frac{\|\vdq Q-\vdq Q_0\|_{F^R}}{\|\vdq Q_0\|_{F^R}}$, 
`$e_x$' and `$e_Q$' are  the error measurements   defined by 
\begin{equation*}
    e_x = \frac{\|\vdq y_r - \vdq x_r\|_{2^R}}{\|\vdq x_r\|_{2^R}} \text{ and }  e_Q = \frac{\|\vdq Q_0-\lambda \vdq u\vdq u^*\|_{F^R}}{\|\vdq Q_0\|_{F^R}}, 
\end{equation*}
where $\vdq x_r = \vdq x\frac{\dq x_i}{|\dq x_i|}$, $i=\mathop{\arg\max}\limits_{i=1,\dots,n} |x_i|$, $\vdq y = \vdq u\sqrt{\lambda}$, and $\vdq y_r = \vdq y\frac{\dq y_i}{|\dq y_i|}$. 
`$n_{\text{iter}}$'  and `time (s)'  are the number of iterations and  the total CPU time  in seconds for implementing the block coordinate descent method in Algorithm \ref{alg:SLAM} respectively. 

\begin{table}
    \centering 
    \caption{Numerical results for a five-point  circle SLAM problem with different noise levels} 
    \begin{tabular}{|c|cccc||c|cccc|}
    \hline
    $l_{\text{noise}}$ & $e_x$ & $e_Q$ & $n_{\text{iter}}$ & time (s) & 
    $l_{\text{noise}}$ & $e_x$ & $e_Q$ & $n_{\text{iter}}$ & time (s) \\\hline
        0\% &   3.85e$-$6 &    4.72e$-$6  & 44.7&   0.1872 &
    5\% &    1.68e$-$2  &   2.20e$-$2 & 144.3&    0.5525\\
    10\% &    3.34e$-$2 &   4.40e$-$2 & 151.6&    0.5958&
    15\% &     4.49e$-$2 &   6.58e$-$2 & 156.1&    0.6489\\
    20\% &    6.63e$-$2  &  8.76e$-$2 & 159.0&    0.6531&
    25\% &     8.25e$-$2 &   1.09e$-$1&  161.7&    0.6949\\
    \hline   
    
    \end{tabular}
    \label{tab:SLAM_5pcircle}
\end{table}
From Table~\ref{tab:SLAM_5pcircle}, we can see that Algorithm~\ref{alg:SLAM} can compute all problems in one minute with a reasonable error. 
Both the error in the vector $\vdq x$ and the matrices $\vdq Q$ are increasing with the noise levels.

We continue to consider the SLAM problem of random graphs. Suppose the number of vertices is $n$  and the observation rate is $s=\frac{m}{n^2}$.  
The numerical results are shown in Table~\ref{tab:SLAM_random}. 
When $n=10$, all results can be obtained in less than one second. 
As the  observation ratio increases, the error in $\vdq x$ and $\vdq Q$ are both decreasing and the number of iterations is decreasing. 
When the observation ratio is above 50\%, then the error in both $\vdq x$ and $\vdq Q$ is near $10^{-6}$. 
When $n=100$, we have similar observations and the error in both $\vdq x$ and $\vdq Q$ is near $10^{-6}$ when the observation ratio is above 30\%.

\begin{table}
    \centering 
    \caption{Numerical results for random graph SLAM problems with  different observation  ratios} 
    \begin{tabular}{|c|cccc||c|cccc|}
    \hline
     \multicolumn{10}{|c|}{$n=10$} \\ 
    \hline 
    $s$ & $e_x$ & $e_Q$ & $n_{\text{iter}}$ & time (s) & 
     $s$ & $e_x$ & $e_Q$ & $n_{\text{iter}}$ & time (s)    \\ 
     \hline
   10\% &   6.62e$-$1 &     4.63e$-$1 &   168.3    &  0.8036 & 
    20\% &    1.50e$-$1 &     1.37e$-$1  &  149.0  &   0.6856 \\
    30\% &    4.29e$-$3  &    4.24e$-$3  &  122.3 &    0.5492 & 
    40\% &     5.74e$-$5  &    5.48e$-$5  &   86.3 &    0.3928 \\
    50\% &     2.35e$-$6 &   1.55e$-$6  &   44.0 &    0.2137 & 
    60\% &    1.28e$-$6  &  1.05e$-$6   &   33.7 &      0.1784 \\ 
    \hline 
     \multicolumn{10}{|c|}{$n=100$} \\ \hline 
    5\% &  3.56e$-$1  &  1.81e$-$1 & 1000 &     14.14 & 
    8\% &  3.95e$-$2  & 3.31e$-$2 & 1000 &     14.58 \\
    10\% &   1.13e$-$2 &    1.15e$-$2 & 987.3 &  14.28 &
    20\% &   8.07e$-$5 &    7.83e$-$5 &  497.0   &   7.171\\
    30\% &    7.78e$-$7  &  8.45e$-$7 &   104.7 &    1.613 &
    40\% &   1.24e$-$7  &  1.40e$-$7   & 39.0 &    0.646\\ 
    \hline  
    \end{tabular}
    \label{tab:SLAM_random}
\end{table}

\bigskip
\section{Final Remarks}\label{sec:finial}

In this paper, we proposed a
power method for computing the dominant eigenvalue of a dual quaternion Hermitian matrix and showed the convergence and convergence rate. 
We used the Laplacian matrices from circles and random graphs to demonstrate the efficiency of the power method. 
Then we further studied the SLAM problem by a two blocks coordinate descent approach, where one block  has a closed form solution and the other block can be obtained by the power method.  
Our results fill the blank of dual quaternion matrix computation and build 
relationship between the dual quaternion matrix theory and the SLAM problem.   Adding the relation between the dual quaternion matrix theory and the formation control problem, established in \cite{QWL22}, we see that the dual quaternion matrix has a solid application background and is worth being further studied.

Several issues remain to be solved.  
\begin{enumerate}
    \item When the multiplicity of the dominant eigenvalue is greater than one, and there exist at least two eigenvalues such that the dual parts are different, how to compute the eigenpairs numerically?
    \item We suspect that the eigenvalues of the Laplacian matrix of a circle are all nonnegative real values, and they are independent of the choice of $\vdq q$. Does this hold?
    \item How to build the convergence analysis of the   block coordinate method for the dual quaternion optimization problem? 
\end{enumerate}
\bigskip



\section*{Data Availability Statement}
\bigskip



 The authors confirm that the data supporting the findings of this study are available within the article and its supplementary materials.
\bigskip



\end{document}